\documentclass{article}
\usepackage{pstricks}

\def\diagram#1{\def\normalbaselines {\lineskip=5pt\baselineskip=0pt
\lineskiplimit=1pt}\matrix{#1}}
\def\hfl#1{\smash{\mathop{\hbox to 10mm{\rightarrowfill}}\limits^{\textstyle
#1}}}
\def\vfl#1{\llap{$\scriptstyle#1$}\left\downarrow\vbox to 6mm{}\right.}

\newtheorem{proposition}[equation]{Proposition}
\newtheorem{corollary}[equation]{Corollary}
\newtheorem{theorem}[equation]{Theorem}
\newtheorem{exa}[equation]{Example}
\newtheorem{ex}[equation]{Exercise}
\newtheorem{s-ex}[equation]{Side-exercise}
\newtheorem{exas}[equation]{Examples}
\newtheorem{lemma}[equation]{Lemma}
\newtheorem{sublemma}[equation]{Sublemma}
\newtheorem{remar}[equation]{Remark}
\newtheorem{remars}[equation]{Remarks}
\newtheorem{nota}[equation]{Notation}
\newtheorem{sremar}[equation]{Side-remark}
\newtheorem{definitio}[equation]{Definition}

\newenvironment{remark}{\begin{remar} \rm }{\end{remar}}

\newenvironment{examples}{\begin{exas} \rm }{\end{exas}}
\newenvironment{example}{\begin{exa} \rm }{\end{exa}}
\newenvironment{definition}{\begin{definitio} \rm }{\end{definitio}}

\newcommand{\CC}{{\bf C}}

\newcommand{\NN}{{\bf N}}
\newcommand{\QQ}{{\bf Q}}
\newcommand{\RR}{{\bf R}}

\newcommand{\ZZ}{{\bf Z}}
\newcommand{\CA}{{\cal A}}
\newcommand{\CB}{{\cal B}}

\newcommand{\CF}{{\cal F}}

\newcommand{\CP}{{\cal P}}

\newcommand{\CO}{{\cal O}}

\newcommand{\bp}{\noindent {\sc Proof: }}
\newcommand{\eop}{\nopagebreak
			\hspace*{\fill}{$\diamond$}
			\medskip}

\newcommand{\fcp}{\begin{pspicture}[.2](0,0)(.5,.4)
\psline{-}(0.05,0)(.45,.4)
\psline{-}(.45,0)(.3,.15)
\psline{-}(.2,.25)(0.05,.4)
\end{pspicture}
}

\newcommand{\fcm}{\begin{pspicture}[.2](0,0)(.5,.4)

\psline{-}(0.45,0)(.05,.4)
\psline{-}(.05,0)(.2,.15)
\psline{-}(.3,.25)(0.45,.4)
\end{pspicture}
}

\newcommand{\fcpoc}{\begin{pspicture}[.2](0,0)(.6,.45)
\psline{->}(0.05,0)(.45,.4)
\psline{-}(.45,0)(.3,.15)
\psline{->}(.2,.25)(0.05,.4)
\rput(.25,0){\tiny C}
\rput(.55,0){\tiny C}
\end{pspicture}
}

\newcommand{\fcmoc}{\begin{pspicture}[.2](0,0)(.6,.45)
\psline{->}(0.45,0)(.05,.4)
\rput(.25,0){\tiny C}
\rput(.55,0){\tiny C}
\psline{-}(.05,0)(.2,.15)
\psline{->}(.3,.25)(0.45,.4)
\end{pspicture}
}

\newcommand{\ass}{\begin{pspicture}[.2](0,0)(.5,.4)
\psline{-}(0.05,0)(.05,.4)
\psline{-}(.15,0)(.35,.4)
\psline{-}(.45,0)(0.45,.4)
\end{pspicture}
}

\newcommand{\assnumi}{\begin{pspicture}[.2](0,0)(.5,.4)
\psline{-}(0.05,0)(.05,.4)
\rput(.05,-.1){\tiny 3}
\psline{-}(.15,0)(.35,.4)
\rput(.15,-.1){\tiny 2}
\psline{-}(.45,0)(0.45,.4)
\rput(.45,-.1){\tiny 1}
\end{pspicture}
}

\newcommand{\snake}{\begin{pspicture}[.2](0,0)(.5,.4)
\pscurve{-}(.05,.4)
(.15,.1)(.35,.3)(.45,0)
\end{pspicture}
}

\newcommand{\bigsnake}{\begin{pspicture}[0.3](0,-.2)(.5,.6)
\pscurve{-}(.05,.6)(.05,0)(.1,-.05)
(.15,0)(.35,.4)(.4,.45)
(0.45,.4)(.45,-.2)
\end{pspicture}
}

\newcommand{\trib}{\begin{pspicture}[.2](0,0)(.8,.4)
\psline{-}(0.05,0)(.05,.4)
\psline{-}(.25,0)(.25,.4)
\rput(.5,.2){$\dots$}
\psline{-}(.75,0)(0.75,.4)
\end{pspicture}
}

\newcommand{\tribr}{\begin{pspicture}[.2](0,0)(1.05,.4)
\psline{-}(0.05,0)(.05,.4)
\psline{-}(.25,0)(.25,.4)
\rput(.5,.2){$\dots$}
\psline{-}(.75,0)(0.75,.4)
\rput(.85,0){\tiny r}
\end{pspicture}
}

\newcommand{\trit}{\begin{pspicture}[.2](0,0)(.5,.4)
\psline{-}(0.05,0)(.05,.4)
\psline{-}(.25,0)(.25,.4)
\psline{-}(.45,0)(0.45,.4)
\end{pspicture}
}

\newcommand{\trid}{\begin{pspicture}[.2](0,0)(.3,.4)
\psline{-}(0.05,0)(.05,.4)
\psline{-}(.25,0)(.25,.4)
\end{pspicture}
}

\newcommand{\itrid}{\begin{pspicture}[.2](0,0)(.1,.25)
\psline{-}(0,0)(0,.25)
\psline{-}(.1,0)(.1,.25)
\end{pspicture}
}

\newcommand{\ptriu}{\begin{pspicture}[.2](0,0)(.1,.25)
\psline{-}(0.05,0)(0.05,.25)
\end{pspicture}
}

\newcommand{\triq}{\begin{pspicture}[.2](0,0)(.4,.4)
\psline{-}(0.05,0)(.05,.4)
\psline{-}(0.15,0)(.15,.4)
\psline{-}(.25,0)(.25,.4)
\psline{-}(0.35,0)(.35,.4)
\end{pspicture}
}

\newcommand{\tridn}{\begin{pspicture}[.2](0,0)(.4,.4)
\psline{->}(0.05,0)(.05,.4)
\rput(.1,0){\tiny 1}
\psline{->}(.25,0)(.25,.4)
\rput(.3,0){\tiny 2}
\end{pspicture}
}

\newcommand{\triu}{\begin{pspicture}[.2](0,0)(.2,.4)
\psline{-}(0.1,0)(.1,.4)
\end{pspicture}
}

\newcommand{\itriu}{\begin{pspicture}[.2](0.1,0)(.2,.25)
\psline{-}(0.1,0)(.1,.25)
\end{pspicture}
}

\newcommand{\DT}{\begin{pspicture}[.2](0,0)(1,1)
\pscurve{-}(.1,.55)(0.05,.6)(0.5,1)(.95,.6)(.95,.4)(0.5,0)(0.05,.4)
(0.25,.6)(0.5,.8)(0.75,.6)(0.75,.4)(0.5,.2)(0.25,.4)(0.2,.45)
\end{pspicture} }

\newcommand{\dcup}{\begin{pspicture}[.2](0,0)(1,.4)
\pscurve{-}(0.05,.4)(0.5,0)(.95,.4)
\pscurve{-}(0.25,.4)(0.5,.2)(0.75,.4)
\end{pspicture} }

\newcommand{\dcap}{\begin{pspicture}[.2](0,0)(1,.4)
\pscurve{-}(0.05,0)(0.5,.4)(.95,0)
\pscurve{-}(0.25,0)(0.5,.2)(0.75,0)
\end{pspicture} }

\newcommand{\dcupn}{\begin{pspicture}[.2](0,0)(1,.4)
\pscurve{->}(0.05,.4)(0.5,0)(.95,.4)
\rput(.15,.35){\tiny 1}
\pscurve{->}(0.25,.4)(0.5,.2)(0.75,.4)
\rput(.4,.35){\tiny 2}
\end{pspicture} }

\newcommand{\dcapn}{\begin{pspicture}[.2](0,0)(1,.4)
\pscurve{->}(0.05,0)(0.5,.4)(.95,0)
\rput(.15,0){\tiny 2}
\pscurve{->}(0.25,0)(0.5,.2)(0.75,0)
\rput(.35,0){\tiny 1}
\end{pspicture} }

\newcommand{\dcapno}{\begin{pspicture}[.2](0,0)(1,.4)
\pscurve{->}(0.05,0)(0.5,.4)(.95,0)
\rput(.15,0){\tiny 1}
\pscurve{->}(0.25,0)(0.5,.2)(0.75,0)
\rput(.35,0){\tiny 2}
\end{pspicture} }

\newcommand{\dcapnof}{\begin{pspicture}[.2](0,0)(1,.4)
\pscurve{->}(0.95,0)(0.5,.4)(.05,0)
\rput(.85,0){\tiny 1}
\pscurve{->}(0.75,0)(0.5,.2)(0.25,0)
\rput(.65,0){\tiny 2}
\end{pspicture} }

\newcommand{\huit}{\begin{pspicture}[.2](0,-.2)(.5,.6)
\pscurve{-}(.3,.25)(0.45,.4)(.25,.6)(.05,.4)(0.45,0)(.25,-.20)(.05,0)(.2,.15) 
\end{pspicture} }

\newcommand{\huito}{\begin{pspicture}[.2](0,-.2)(.5,.6)
\pscurve{-}(.3,.15)(0.45,0)(.25,-.2)(.05,0)(0.45,.4)(.25,.6)(.05,.4)(.2,.25) 
\end{pspicture} }

\newcommand{\scirc}{\begin{pspicture}[.2](0,0)(.4,.4)
\pscircle(.2,.2){.15}
\end{pspicture} }

\newcommand{\dmer}{\begin{pspicture}[0.3](0,0)(.5,1.6)
\psline{->}(.4,0)(.4,1.5)
\psline[linestyle=dashed,dash=2pt 1pt]{*-*}(.1,.8)(.4,.8)
\psline[linestyle=dashed,dash=2pt 1pt]{-*}(.1,.8)(.4,.4)
\psline[linestyle=dashed,dash=2pt 1pt]{-*}(.1,.8)(.4,1.1)
\end{pspicture}}

\newcommand{\firstleg}{\begin{pspicture}[.2](0,0)(.35,.3)
\psline{-}(0.1,0)(.1,.3)
\psline[linestyle=dashed,dash=2pt 1pt]{-*}(.25,.1)(0.1,.1)
\end{pspicture}
}

\newcommand{\secondleg}{\begin{pspicture}[.2](0,0)(.8,.3)
\pscurve[linestyle=dashed,dash=2pt 1pt]{-*}(.1,.1)(.3,.15)(.7,.2)(.55,.1)
\psline{-}(.55,0)(0.55,.3)
\end{pspicture}
}

\newcommand{\tata}{\begin{pspicture}[0.2](0,0)(.5,.4)
\pscircle(0.25,0.2){.2}
\psline[linestyle=dashed,dash=2pt 1pt]{*-*}(0.05,.2)(.45,.2)
\end{pspicture}}

\newcommand{\cirt}{\begin{pspicture}[0.4](-.9,-.8)(.9,.8)
\psarc{->}(0,0){.8}{180}{20}
\psarc{-}(0,0){.8}{20}{180}
\SpecialCoor
\psline[linestyle=dashed,dash=2pt 1pt]{*-*}(.8;135)(0,.2)(.8;45)
\psline[linestyle=dashed,dash=2pt 1pt]{*-*}(.8;-135)(0,-.2)(.8;-45)
\psline[linestyle=dashed,dash=2pt 1pt]{*-*}(0,-.2)(0,.2)
\end{pspicture}}

\newcommand{\cirtd}{\begin{pspicture}[0.4](-.9,-.8)(.9,.8)
\psarc{->}(0,0){.8}{180}{20}
\psarc{-}(0,0){.8}{20}{180}
\SpecialCoor
\psline[linestyle=dashed,dash=2pt 1pt]{*-*}(.8;135)(0,.2)(.8;45)
\psline[linestyle=dashed,dash=2pt 1pt]{*-*}(.8;-135)(0,-.2)(.8;-45)
\psline[linestyle=dashed,dash=2pt 1pt]{*-*}(0,-.2)(0,.2)
\psline[linestyle=dashed,dash=2pt 1pt]{*-*}(.3,0)(.8;30)
\psline[linestyle=dashed,dash=2pt 1pt]{-*}(.3,0)(.8;0)
\psline[linestyle=dashed,dash=2pt 1pt]{-*}(.3,0)(.8;-30)
\end{pspicture}}

\newcommand{\cirtdp}{\begin{pspicture}[0.4](-.9,-.8)(.9,.8)
\psarc{->}(0,0){.8}{180}{20}
\psarc{-}(0,0){.8}{20}{180}
\SpecialCoor
\psline[linestyle=dashed,dash=2pt 1pt]{*-*}(.8;135)(0,.2)(.8;45)
\psline[linestyle=dashed,dash=2pt 1pt]{*-*}(.8;-135)(0,-.2)(.8;-45)
\psline[linestyle=dashed,dash=2pt 1pt]{*-*}(0,-.2)(0,.2)
\psline[linestyle=dashed,dash=2pt 1pt]{*-*}(0,.4)(.8;60)
\psline[linestyle=dashed,dash=2pt 1pt]{-*}(0,.4)(.8;90)
\psline[linestyle=dashed,dash=2pt 1pt]{-*}(0,.4)(.8;120)
\end{pspicture}}

\newcommand{\smin}{\begin{pspicture}[.2](0,0)(.5,.4)
\pscurve{-}(.05,.4)(0.25,0)(.45,.4)
\end{pspicture}
}

\newcommand{\sminc}{\begin{pspicture}[.2](0,-.2)(.5,.4)
\psline{-}(.3,.25)(.45,.4)
\pscurve{-}(.05,.4)(.25,.2)(0.45,0)(.25,-.20)(.05,0)(.2,.15) 
\end{pspicture} }

\newcommand{\smax}{\begin{pspicture}[.2](0,0)(.5,.4)
\pscurve{-}(.05,0)(0.25,0.4)(.45,0)
\end{pspicture}
}

\newcommand{\ismax}{\begin{pspicture}[.2](0,0)(.2,.2)
\pscurve{-}(0,0)(0.1,0.2)(.2,0)
\end{pspicture}
}

\newcommand{\ssmaxn}{\begin{pspicture}[.2](0,0)(.4,.4)
\pscurve{-}(.25,.2)(0.35,.3)(.2,.4)(.05,.3)(.2,.15)(0.35,0)
\psline{-}(.05,0)(.15,.1) 
\end{pspicture} }

\newcommand{\smaxn}{\begin{pspicture}[.2](0,0)(.5,.6)
\pscurve{-}(.3,.25)(0.45,.4)(.25,.6)(.05,.4)(.25,.2)(0.45,0)
\psline{-}(.05,0)(.2,.15) 
\end{pspicture} }

\newcommand{\hbelt}{\begin{pspicture}[.2](-.05,-.1)(.5,.5)
\pscurve{-}(0.05,.5)(.1,0)(.2,.1)
\pscurve{-}(.3,.2)(.35,.25)(.25,.35)(.15,.25)(.25,.15)(0.45,-.1) 
\end{pspicture} }

\newcommand{\hbbelt}{\begin{pspicture}[.3](-.05,-.1)(.5,.7)

\pscurve{-}(.3,.4)(.35,.45)(0.25,.55)(.15,.45)(.25,.35)(.35,.25)(.3,.2)
\pscurve{-}(0.05,.7)(.1,0)(.2,.1)
\pscurve{-}(.2,.3)(.15,.25)(.25,.15)(0.45,-.1) 
\end{pspicture} }

\newcommand{\smaxo}{\begin{pspicture}[.2](0,0)(.5,.4)
\pscurve{->}(.05,0)(0.25,0.4)(.45,0)
\end{pspicture}
}

\newcommand{\smaxno}{\begin{pspicture}[.2](0,0)(.5,.6)

\pscurve{->}(.3,.25)(0.45,.4)(.25,.6)(.05,.4)(.25,.2)(0.45,0)
\psline{-}(.05,0)(.2,.15) 
\end{pspicture} }

\title{About the uniqueness and the denominators of the Kontsevich integral}
\author{Christine Lescop \thanks{CNRS, Institut Fourier (UMR 5582)}}

\begin{document}

\maketitle
\begin{abstract}
We refine a Le and Murakami uniqueness theorem for the Kontsevich Integral in order to specify the relationship between the two (possibly equal) main universal Vassiliev link invariants: the Kontsevich Integral and the perturbative expression of the Chern-Simons theory.

As a corollary, we prove that the Altschuler and Freidel anomaly $\alpha \in \CA(S^1)$  -that groups the Bott and Taubes anomalous terms- is a combination of diagrams with two univalent vertices and we explicitly define the isomorphism of $\CA$ which 
transforms the Kontsevich integral into the Poirier limit
of the perturbative expression of the Chern-Simons theory for framed links, as a function of $\alpha$. 

As a consequence of this corollary, we use the Poirier estimates on the denominators of the perturbative expression of the Chern-Simons theory to show that the denominators of the degree n part of the Kontsevich integral of framed
links divide into $(2!3! \dots (n-5)!)(n-5)!3^2(3n-4)!2^{2n+2}$ for $n \geq 5$.
\end{abstract}

\begin{flushleft}
\begin{small}
A.M.S. subject classification: 57M27, 57M25, 17B37.
\end{small}
\end{flushleft}

\section{Introduction}
There are essentially two universal Vassiliev invariants of links,
the {\em Kontsevich integral,\/} and the {\em perturbative expression of the Chern-Simons theory\/} studied by Guadagnini Martellini and Mintchev, Bar-Natan, Kontsevich, Bott and Taubes~\cite{bt}, D. Thurston~\cite{th}, Altschuler and Freidel~\cite{af}, Yang, Poirier~\cite{p1,p2,p3}... 
The question of knowing whether the two invariants
coincide or not is still open to mathematicians despite the substantial progress of Sylvain Poirier who in particular reduced this question to the computation of the {\em anomaly} of Bott, Taubes, Altschuler and Freidel which is an element
of the space of Feynman diagrams $\CA(S^1)$.

In this article, after specifying all our definitions and notations, we first review all
the main properties that are shared by the Kontsevich integral and the Poirier limit of the perturbative expression of the Chern-Simons theory. We call an invariant of framed links satisfying all these properties {\em a good monoidal functor from the category
of framed q-tangles to $\CA$.\/} Second, we refine a theorem of
Le and Murakami \cite[Theorem 8]{lm} inspired by Drinfeld and Kontsevich that could be
stated as follows: {\em A good monoidal functor that varies like the Kontsevich integral $Z^K$ under a framing change must coincide with $Z^K$ on framed links.\/} Our refinement is stated in Theorem~\ref{mainth}. Roughly speaking, we list the possible variations of good monoidal functors under framing changes and we define some special type of isomorphisms of $\CA$, {\em the $\Psi(\beta)$}, so that, when restricted to framed links, any good monoidal functor is of the form $\Psi(\beta) \circ Z^K$.

As a corollary, we explicitly define the isomorphism of $\CA$ which 
transforms the Kontsevich integral into the Poirier limit \cite{p3}
of the perturbative expression of the Chern-Simons theory, as a function of the Bott and Taubes anomaly. This corollary allows us to use the Poirier estimates on the denominators of the perturbative expression of the Chern-Simons theory to show that the denominators of the degree n part of the Kontsevich integral of framed links divide into $(2!3! \dots (n-5)!)(n-5)!3^2(3n-4)!2^{2n+2}$ for $n \geq 5$. These denominators have more 
prime factors than the denominators $(2!3! \dots n!)^4(n+1)!$ of Le \cite{le} but they are smaller.

I thank Christophe Champetier, Gregor Masbaum,  and especially Dylan Thurston and Pierre Vogel for useful conversations.


\begin{definition}
\label{defbetas}
Here, we use standard notation related to spaces of Feynman diagrams that will be recalled in Section \ref{secdia}. We say that an element $\beta=(\beta_n)_{n \in {\bf N}}$ in $\CA(S^1)$ is a {\em two-leg element\/} if, for any $n \in {\bf N}$, $\beta_n$ is a combination of diagrams with two univalent 
vertices.
 
Forgetting $S^1$ from such a two-leg element gives rise to a unique series $\beta^s$ of diagrams
with two distinguished univalent vertices $v_1$ and $v_2$, such that $\beta^s$ is symmetric with respect to the exchange of $v_1$ and $v_2$. (According to
\cite[Corollary 4.2]{vo}, all two-leg elements are symmetric with respect to this symmetry modulo the standard AS and IHX relations. This is reproved as 
Lemma~\ref{lemsymvo} below.)

If $\Gamma$ is a chord diagram, then $\Psi(\beta)(\Gamma)$ is defined by replacing each chord
by $\beta^s$. As it will be seen after Lemma~\ref{lemins}, $\Psi(\beta)$ is a well-defined morphism of topological vector spaces from $\CA(M)$ to $\CA(M)$ for any one-manifold $M$, and  $\Psi(\beta)$ is an isomorphism as soon as $\beta_1 \neq 0$.
\end{definition}
 
All the natural
definitions used in the following definition will be specified in Section \ref{secfunc}.

\begin{definition}
A {\em good monoidal functor\/} is 
 a monoidal functor $Z$ from the category of framed q-tangles (up to framed isotopy) to the category
of spaces of Feynman diagrams which satisfies:
\begin{enumerate}
\item For any q-tangle $T$, the degree zero part $Z_0(T)$ of $Z(T)$ is 1.
\item $Z$ is compatible with the deletion of a component.
\item $Z$ is compatible with the duplication of a {\em regular\/} component, that is a component that can be represented without horizontal tangent vectors.
\item $Z$ is invariant under the 180-degree rotation around the vertical axis. \item Let $s_h$ be the orthogonal symmetry with respect to the horizontal plane and let  $s_v$ be the orthogonal symmetry with respect to the  blackboard plane.
Let $\sigma_v$ and $\sigma_h$ be the two endomorphisms of the topological vector spaces $\CA(S^1)$ such that $\sigma_v(z[\Gamma]) = (-1)^d \overline{z}[\Gamma]$
and $\sigma_h(z[\Gamma]) = (-1)^d z[\Gamma]$, where
$z$ is a complex number and $[\Gamma]$ is the image of a degree $d$ diagram $\Gamma$ in $\CA(S^1)$. Then, for any framed knot $K(S^1)$, 
$$Z \circ s_v(K)=\sigma_v \circ Z(K)$$ and $$Z \circ s_h(K) =\sigma_h \circ Z(K)$$ 
\item The element $a^Z \in \CA(S^1)$ such that $a^Z_0=0$ and $$Z\left(\smaxno \right) = \exp(a^Z)Z\left(\smaxo \right)$$ has a nonzero degree one part $$a^Z_1 \neq 0.$$
\end{enumerate}
\end{definition}

\noindent This article is devoted to proving the following result.
\begin{theorem}
\label{mainth}
If $Z$ is a good monoidal functor as above, then $a^Z$ is a
 two-leg element of $\CA(S^1)$, such that for any integer $i$, $a^Z_{2i}=0$ and $a^Z_{2i+1}$ is a combination of 
diagrams with real coefficients, and,
for any framed link $L$, 
$$Z(L) = \Psi(2a^Z)(Z^K(L))$$ where $Z^K$ denotes the Kontsevich integral of framed links (denoted by $\hat{Z}_f$ in \cite{lm} and by $Z$ in \cite{les}).
\end{theorem}

It follows from \cite{lm} and \cite{p3} that both the Kontsevich integral $Z^K$ and the Poirier limit integral $Z^{\ell}$ satisfy
the hypotheses of our theorem with $a^{Z^K}=\frac12\tata$ and $a^{Z^{\ell}}=\alpha$. Thus, we can state the following corollary using the notation of \cite{p3}.

\begin{corollary}
The anomaly $\alpha$ is a two-leg element of $\CA(S^1)$.
For any framed link $L$, the Poirier limit integral $Z^{\ell}(L)$ is equal to
$\Psi(2\alpha)(Z^K(L))$.
\end{corollary}

\begin{remark}
The Poirier limit integral and the perturbative expression of the Chern-Simons theory coincide for the framed links whose components all have both zero framing and zero Gauss integral.
\end{remark}

As it will be shown in Section~\ref{secden}, the denominator estimates of Sylvain Poirier for the perturbative expression of the Chern-Simons theory induce the following estimates for the denominators of
the Kontsevich integral.

\begin{corollary}
\label{corcoco}
Let $L$ be a framed link and let $Z^{K}_n(L)$ denote the degree $n$ part of its Kontsevich integral. Let
$$d(n)=\left\{\begin{array}{ll}
(3n-4)!2^{3n-4} & \mbox{if}\;3 \leq n \leq 8\\
(2!3! \dots (n-5)!)(n-5)!3^2(3n-4)!2^{2n+1} &  \mbox{if}\;n \geq 9 \;\mbox{and n is odd}\\
(2!3! \dots (n-6)!)(n-6)!3^2(3n-4)!2^{2n-1} &  \mbox{if}\;n \geq 10 \;\mbox{and n is even}
\end{array} \right. $$
Then for all $n \geq 3$, $d(n)Z^{K}_n(L)$ is an integral combination of chord diagrams.
\end{corollary}

\section{First definitions on the spaces of Feynman diagrams}
\label{secdia}
\setcounter{equation}{0}

\begin{definition}
\label{defdia}
Let $M$ be a one-manifold and let $X$ be a finite set.
A {\em diagram\/} $\Gamma$ with support $M \cup X$ is a finite uni-trivalent graph $\Gamma$ such that
every connected component of $\Gamma$ has at least one univalent vertex, 
equipped with:
\begin{enumerate}
\item a partition of the set $U$ of univalent vertices of $\Gamma$ also called {\em legs} of $\Gamma$ into two (possibly empty) subsets $U_X$ and $U_M$,
\item a function $f$ from $U_X$ to $X$,
\item an isotopy class of injections $i$ of $U_M$ into the interior of $M$,
\item an {\em orientation\/} of every trivalent vertex, that is a cyclic order
on the set of the three half-edges which meet at this vertex,
\item an {\em orientation\/} of every univalent vertex $u$ of $U_M$, that is
a cyclic order on the set made of the half-edge that contains $u$ and the two sides of $i(u)$ on $M$. 
\end{enumerate}
\end{definition}

Such a diagram $\Gamma$ is represented by a planar immersion of $\Gamma \cup M$ where the univalent vertices of $U_M$ are located at their images under $i$, the one-manifold $M$ is represented by solid lines, whereas
the diagram $\Gamma$ is dashed. The vertices are represented by big points. The local orientation of a 
vertex is represented by the counterclockwise order of the three half-edges
(solid or dashed) that meet at it.

Here is an example of a diagram $\Gamma$ on the disjoint 
union $M=S^1 \coprod S^1$ of two circles:

$$
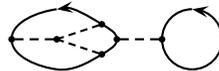
\begin{pspicture}[.2](0,0)(3,1)
\psecurve{->}(1.3,.7)(.7,.9)(.1,.5)(.7,.1)(1.3,.3)(1.5,.5)(1.3,.7)(.7,.9)(.1,.5)
\psarc{->}(2.5,.5){.4}{-90}{90}
\psarc{-}(2.5,.5){.4}{90}{-90}
\psline[linestyle=dashed,dash=3pt 2pt]{*-*}(1.5,.5)(2.1,.5)
\psline[linestyle=dashed,dash=3pt 2pt]{*-*}(.1,.5)(.7,.5)
\psline[linestyle=dashed,dash=3pt 2pt]{-*}(.7,.5)(1.3,.3)
\psline[linestyle=dashed,dash=3pt 2pt]{-*}(.7,.5)(1.3,.7)
\end{pspicture}$$

The {\em degree} of such a diagram is 
half the number of all its vertices. 

A {\em chord diagram\/} is a diagram on a one-manifold $M$ ($X=\emptyset$) without trivalent vertices.

Let $\CA_n(M \cup X)$ denote the quotient of the complex vector space generated by the degree n diagrams
on $M \cup X$ by the following relations AS, STU and IHX:

$$ {\rm AS:}  \begin{pspicture}[.2](0,-.2)(1,1)
\psline[linestyle=dashed,dash=3pt 2pt]{*-}(.5,.5)(.5,0)
\psline[linestyle=dashed,dash=3pt 2pt]{-}(.1,.9)(.5,.5)
\psline[linestyle=dashed,dash=3pt 2pt]{-}(.9,.9)(.5,.5)
\end{pspicture}
+
\begin{pspicture}[.2](0,-.2)(1,1)
\psline[linestyle=dashed,dash=3pt 2pt]{*-}(.5,.5)(.5,0)
\pscurve[linestyle=dashed,dash=3pt 2pt]{-}(.1,.9)(.7,.7)(.5,.5)
\pscurve[linestyle=dashed,dash=3pt 2pt]{-}(.9,.9)(.3,.7)(.5,.5)
\end{pspicture}
=0 \;\;{\rm and} \;\;
\begin{pspicture}[.2](0,-.2)(.8,1)
\psline[linestyle=dashed,dash=3pt 2pt]{-*}(.4,.9)(.4,.1)
\psline{-}(.1,.1)(.7,.1)
\end{pspicture}
+
\begin{pspicture}[.2](0,-.2)(.8,1)
\pscurve[linestyle=dashed,dash=3pt 2pt]{-*}(.4,.9)(.3,-.1)(.45,.1)
\psline{-}(.1,.1)(.7,.1)
\end{pspicture}
=0$$

$${\rm STU:} \begin{pspicture}[.2](0,-.2)(1,1)
\psline[linestyle=dashed,dash=3pt 2pt]{*-*}(.5,.5)(.5,.1)
\psline[linestyle=dashed,dash=3pt 2pt]{-}(.3,.9)(.5,.5)
\psline[linestyle=dashed,dash=3pt 2pt]{-}(.7,.9)(.5,.5)
\psline{-}(.1,.1)(1,.1)
\end{pspicture}
=
\begin{pspicture}[.2](0,-.2)(1,1)
\psline[linestyle=dashed,dash=3pt 2pt]{-*}(.3,.9)(.3,.1)
\psline[linestyle=dashed,dash=3pt 2pt]{-*}(.7,.9)(.7,.1)
\psline{-}(.1,.1)(1,.1)
\end{pspicture}
-
\begin{pspicture}[.2](0,-.2)(1,1)
\psline[linestyle=dashed,dash=3pt 2pt]{-*}(.3,.9)(.7,.1)
\psline[linestyle=dashed,dash=3pt 2pt]{-*}(.7,.9)(.3,.1)
\psline{-}(.1,.1)(1,.1)
\end{pspicture} 
\;\;{\rm equivalent\; under\; AS \;to} \;\;
\begin{pspicture}[.2](0,-.2)(1,1)
\psline[linestyle=dashed,dash=3pt 2pt]{-*}(.2,.9)(.3,.1)
\psline[linestyle=dashed,dash=3pt 2pt]{-*}(.5,.9)(.5,.1)
\psline{-}(.1,.1)(1,.1)
\end{pspicture}
+
\begin{pspicture}[.2](0,-.2)(1,1)
\pscurve[linestyle=dashed,dash=3pt 2pt]{-*}(.2,.9)(.5,-.1)(.7,.1)
\psline[linestyle=dashed,dash=3pt 2pt]{-*}(.5,.9)(.5,.1)
\psline{-}(.1,.1)(1,.1) 
\end{pspicture}
+
\begin{pspicture}[.2](0,-.2)(1,1)
\pscurve[linestyle=dashed,dash=3pt 2pt]{-*}(.2,.9)(.4,.3)(.7,.45)(.5,.6)
\psline[linestyle=dashed,dash=3pt 2pt]{-*}(.5,.9)(.5,.1)
\psline{-}(.1,.1)(1,.1)
\end{pspicture}
=0
$$

$${\rm IHX:} 
\begin{pspicture}[.2](0,-.2)(1,1)
\psline[linestyle=dashed,dash=3pt 2pt]{-*}(.1,1)(.35,.2)
\psline[linestyle=dashed,dash=3pt 2pt]{*-}(.5,.5)(.5,1)
\psline[linestyle=dashed,dash=3pt 2pt]{-}(.75,0)(.5,.5)
\psline[linestyle=dashed,dash=3pt 2pt]{-}(.25,0)(.5,.5)
\end{pspicture}
+
\begin{pspicture}[.2](0,-.2)(1,1)
\pscurve[linestyle=dashed,dash=3pt 2pt]{-*}(.1,1)(.3,.3)(.65,.2)
\psline[linestyle=dashed,dash=3pt 2pt]{*-}(.5,.5)(.5,1)
\psline[linestyle=dashed,dash=3pt 2pt]{-}(.75,0)(.5,.5)
\psline[linestyle=dashed,dash=3pt 2pt]{-}(.25,0)(.5,.5)
\end{pspicture}
+
\begin{pspicture}[.2](0,-.2)(1,1)
\pscurve[linestyle=dashed,dash=3pt 2pt]{-*}(.1,1)(.2,.7)(.7,.7)(.5,.8)
\psline[linestyle=dashed,dash=3pt 2pt]{*-}(.5,.5)(.5,1)
\psline[linestyle=dashed,dash=3pt 2pt]{-}(.75,0)(.5,.5)
\psline[linestyle=dashed,dash=3pt 2pt]{-}(.25,0)(.5,.5)
\end{pspicture}
=0$$

Each of these relations relate diagrams which are identical outside the pictures
where they are like in the pictures. For example, 
AS identifies the sum of two diagrams which only differ by the orientation
at one vertex to zero.

Let $$\CA(M \cup X) =\prod_{n \in {\bf N}}\CA_n(M \cup X)$$ denote the product of the $\CA_n(M \cup X)$ as a topological vector space.
$\CA_0(M \cup X)$ is equal to $\CC$ generated by the empty diagram.

\begin{lemma}
\label{lemcom}
Let $\Gamma_1$ be a diagram  with support $M \cup X$. Assume that $\Gamma_1 \cup M$ is immersed in the plane so that  $\Gamma_1 \cup M$ meets an open annulus $A$ embedded in the plane
exactly along $n+1$ arcs $\alpha_1$, $\alpha_2$, \dots, $\alpha_n$ and $\beta$,
and one vertex $v$ so that:
\begin{enumerate}
\item The $\alpha_i$ may be dashed or solid, they run from a boundary component of $A$ to the other one,
\item $\beta$ is a dashed arc which runs from the boundary of $A$ to
$v \in \alpha_1$,
\item The bounded component $D$ of the complement of $A$ does  contain neither a boundary point of $M$ nor a univalent vertex in $U_X$. 
\end{enumerate}
Let $\Gamma_i$ be the diagram obtained from $\Gamma_1$ by attaching the 
endpoint $v$ of $\beta$ to $\alpha_i$ instead of $\alpha_1$ on the same side,
where
the side of an arc is its side when going from the outside boundary component
of $A$ to the inside one $\partial D$.
Then, we have in $\CA(M \cup X)$
$$\sum_{i=1}^n\Gamma_i=0$$
\end{lemma}
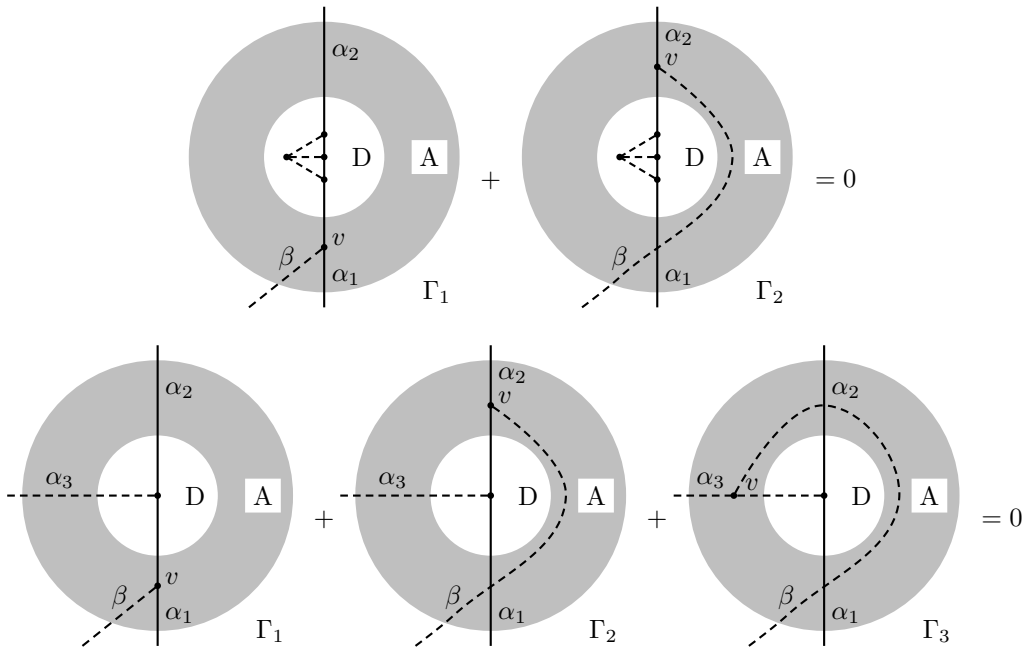
\begin{examples}
$$\begin{pspicture}[.4](0,0)(4,4)
\pscircle[dimen=inner,linewidth=1,linecolor=lightgray](2,2){.8}
\psline[linestyle=dashed,dash=3pt 2pt]{*-*}(1.5,2)(2,2)
\psline[linestyle=dashed,dash=3pt 2pt]{-*}(1.5,2)(2,1.7)
\psline[linestyle=dashed,dash=3pt 2pt]{-*}(1.5,2)(2,2.3)
\psline{-}(2,0)(2,4)
\psline[linestyle=dashed,dash=3pt 2pt]{-*}(1,0)(2,.8)
\rput[l](2.1,.9){$v$}
\rput[l](2.1,.4){$\alpha_1$}
\rput[l](2.1,3.4){$\alpha_2$}
\rput(1.5,.65){$\beta$}
\rput(3.4,2){\psframebox*{A}}
\rput(2.5,2){D}
\rput(3.5,.2){$\Gamma_1$}
\end{pspicture}
+ 
\begin{pspicture}[.4](0,0)(4,4)
\pscircle[dimen=inner,linewidth=1,linecolor=lightgray](2,2){.8}
\psline[linestyle=dashed,dash=3pt 2pt]{*-*}(1.5,2)(2,2)
\psline[linestyle=dashed,dash=3pt 2pt]{-*}(1.5,2)(2,1.7)
\psline[linestyle=dashed,dash=3pt 2pt]{-*}(1.5,2)(2,2.3)
\psline{-}(2,0)(2,4)
\pscurve[linestyle=dashed,dash=3pt 2pt]{-*}(1,0)(1.5,.4)(3,2)(2,3.2)
\rput[l](2.1,3.3){$v$}
\rput[l](2.1,.4){$\alpha_1$}
\rput[l](2.1,3.6){$\alpha_2$}
\rput(1.5,.65){$\beta$}
\rput(3.4,2){\psframebox*{A}}
\rput(2.5,2){D}
\rput(3.5,.2){$\Gamma_2$}
\end{pspicture}
=0
$$

$$
\begin{pspicture}[.4](0,0)(4,4)
\pscircle[dimen=inner,linewidth=1,linecolor=lightgray](2,2){.8}
\psline[linestyle=dashed,dash=3pt 2pt]{-*}(0,2)(2,2)
\rput(.7,2.2){$\alpha_3$}
\psline{-}(2,0)(2,4)
\psline[linestyle=dashed,dash=3pt 2pt]{-*}(1,0)(2,.8)
\rput[l](2.1,.9){$v$}
\rput[l](2.1,.4){$\alpha_1$}
\rput[l](2.1,3.4){$\alpha_2$}
\rput(1.5,.65){$\beta$}
\rput(3.4,2){\psframebox*{A}}
\rput(2.5,2){D}
\rput(3.5,.2){$\Gamma_1$}
\end{pspicture}
+
\begin{pspicture}[.4](0,0)(4,4)
\pscircle[dimen=inner,linewidth=1,linecolor=lightgray](2,2){.8}
\psline[linestyle=dashed,dash=3pt 2pt]{-*}(0,2)(2,2)
\rput(.7,2.2){$\alpha_3$}
\psline{-}(2,0)(2,4)
\pscurve[linestyle=dashed,dash=3pt 2pt]{-*}(1,0)(1.5,.4)(3,2)(2,3.2)
\rput[l](2.1,3.3){$v$}
\rput[l](2.1,.4){$\alpha_1$}
\rput[l](2.1,3.6){$\alpha_2$}
\rput(1.5,.65){$\beta$}
\rput(3.4,2){\psframebox*{A}}
\rput(2.5,2){D}
\rput(3.5,.2){$\Gamma_2$}
\end{pspicture}
+
\begin{pspicture}[.4](0,0)(4,4)
\pscircle[dimen=inner,linewidth=1,linecolor=lightgray](2,2){.8}
\psline[linestyle=dashed,dash=3pt 2pt]{-*}(0,2)(2,2)
\rput(.5,2.2){$\alpha_3$}
\psline{-}(2,0)(2,4)
\pscurve[linestyle=dashed,dash=3pt 2pt]{-*}(1,0)(1.5,.4)(3,2)(2,3.2)(.8,2)
\rput[l](.95,2.15){$v$}
\rput[l](2.1,.4){$\alpha_1$}
\rput[l](2.1,3.4){$\alpha_2$}
\rput(1.5,.65){$\beta$}
\rput(3.4,2){\psframebox*{A}}
\rput(2.5,2){D}
\rput(3.5,.2){$\Gamma_3$}
\end{pspicture}
=0
$$

\end{examples}
\bp
The second example shows that $STU$ is equivalent to this relation when the bounded component $D$ of $\RR^2 \setminus A$ intersects $\Gamma_1$ in the neighborhood of a univalent vertex on $M$. Similarly, $IHX$ is easily seen as given by this relation  when
 $D$ intersects $\Gamma_1$ in the neighborhood of a trivalent vertex. Also note that
$AS$  corresponds to the case when $D$  intersects $\Gamma_1$ along a
dashed or solid arc.
Now for the Bar-Natan \cite[Lemma 3.1]{bn} proof. See also \cite[Lemma 3.3]{vo}. Assume without loss that $v$ is always attached on the right-hand-side of 
the $\alpha$'s. Add to the sum the trivial (by IHX and STU) contribution of the sum of the diagrams
obtained from $\Gamma_1$ by attaching $v$ to each of the three (dashed or solid) half-edges of each vertex $w$
of $\Gamma_1 \cup M$ in $D$ on the left-hand side when the half-edges are oriented towards $w$. Now, group the terms of the obtained sum by edges of $\Gamma_1 \cup M$
where $v$ is attached, and observe that the sum is zero edge by edge by $AS$.
\eop

Assume that a one-manifold $M$ is decomposed as a union of two one-manifolds $M = M_1 \cup M_2$ whose interiors in $M$ do not intersect. Then, we define the
{\em product associated to this decomposition\/}:
$$\CA(M_1) \times \CA(M_2) \longrightarrow \CA(M)$$
as the continuous bilinear map which maps $([\Gamma_1],[\Gamma_2])$ to $[\Gamma_1 \coprod \Gamma_2]$, if $\Gamma_1$ is a diagram with support $M_1$ and if $\Gamma_2$ is a diagram with support $M_2$, where $\Gamma_1 \coprod \Gamma_2$ denotes their disjoint union. In the particular case when $M_1$ and $M_2$ are 
disjoint, this product is sometimes denoted by $\otimes$.

Let $I=[0,1]$ be the compact oriented interval.
Another particular case is the case when $M$ is an ordered union of $p$ intervals
which are seen as vertical $$M=\trib=\{1,2,\dots,p\} \times [0,1]$$
Then if we naturally identify $M$ to $M_1=\{1,2,\dots,p\} \times [0,1/2]$
and to $M_2=\{1,2,\dots,p\} \times [1/2,1]$.
The above process turns $\CA(M)$ into an algebra where the elements with
degree zero part 1 admit an inverse and a unique square root whose degree 0 part is 1.

With each choice of a connected component $C$ of $M$ and of an orientation of $C$, we associate an {\em $\CA(I)$-module structure on $\CA(M)$\/}, that is given by the continuous bilinear map:
$$\CA(I) \times \CA(M) \longrightarrow \CA(M)$$
such that:
If  $\Gamma^{\prime}$ is a diagram with support $M$ and if $\Gamma$ is a diagram with support $I$, then $([\Gamma],[\Gamma^{\prime}])$ is mapped to the class of the diagram obtained by inserting $\Gamma$ along $C$ outside the vertices of $\Gamma$, according to the given orientation. For example, $$\dmer \cirt = \cirtd =\cirtdp $$
As shown in the first example that illustrates  Lemma \ref{lemcom}, the independence of the choice of the insertion locus is a consequence of Lemma~\ref{lemcom} where $\Gamma_1$ is the disjoint union $\Gamma \coprod \Gamma^{\prime}$ and intersects $D$ along $\Gamma \cup I$.
This also proves that $\CA(I)$ is a commutative algebra.
Since the morphism from $\CA(I)$ to $\CA(S^1)$ induced by the identification of the two endpoints of $I$ amounts to quotient out $\CA(I)$ by the relation that identifies two diagrams that are obtained from one another by moving the nearest univalent vertex to an endpoint of I near the other endpoint, a similar application of Lemma \ref{lemcom} also proves that this morphism is an isomorphism from $\CA(I)$ to $\CA(S^1)$. (In this application,
$\beta$ comes from the inside boundary of the annulus.) This identification between
$\CA(I)$ and $\CA(S^1)$ will be used several times.

Let $C$ be a component of $M \cup X$, that is either a point of $X$ or a
connected component of $M$.
Then we define $(r \times C)(M \cup X)$ from $M \cup X$ by replacing $C$
by $r$ copies $C_1, C_2, \cdots, C_r$ of $C$:
$$(r \times C)(M \cup X)= \left((M \cup X) \setminus C \right) \cup \bigcup_{i=1}^r C_i$$

Let $\Gamma$ be a diagram with support $M \cup X$ as in Definition~\ref{defdia}. Let $U_C$ denote
the preimage of $C$ either under $f$ if $C \in X$, or under $i$ if 
$C \subset M$. 
Let $\CP$ be the set of ordered partitions $P = (U_1, \cdots, U_r)$ of $U_C$, that are ordered collections of disjoint $U_j$, with $\emptyset \subseteq U_j
\subseteq U_C$ such that $\bigcup_{j=1}^rU_j=U_C$.
To every $P \in \CP$, we associate the diagram $\Gamma_P$ on  $(r \times C)(M \cup X)$, obtained from $\Gamma$ by changing 
$f$ and $i$ into $\tilde{f}$ and $\tilde{i}$ so that $\tilde{f}=f$ and $\tilde{i}=i$ outside $U_C$ and $\tilde{f}(U_j)=\{C_j\}$ if $C\in X$ and
$\tilde{i}_{|U_j}=\iota_j \circ i_{|U_j}$, where $\iota_j$ is the identification
morphism from $C$ to $C_j$ which also carries the local orientations of the vertices of $U_j$.

The {\em duplication map\/} $(r \times C)_{\ast}$ from $\CA(M \cup X)$ to $\CA((r \times C)(M \cup X))$ is the (well-defined!) morphism of topological vector spaces which maps
 $[\Gamma]$ to $$(r \times C)_{\ast}([\Gamma])=\sum_{P \in \CP} [\Gamma_P]$$

Locally, this reads:
$$\begin{pspicture}[.2](0,0)(.6,.6)
\psline{-}(0.4,0.1)(.4,.5)
\psline[linestyle=dashed,dash=3pt 2pt]{-*}(.1,.3)(.4,.3)
\end{pspicture}
\stackrel{(r \times \ptriu)_{\ast}}{\mapsto}
\begin{pspicture}[.2](-.2,0)(1.4,.6)
\psline{-}(0.3,0.1)(.3,.5)
\psline{-}(0.5,0.1)(.5,.5)
\psline{-}(1.2,0.1)(1.2,.5)
\psline[linestyle=dashed,dash=3pt 2pt]{-*}(.1,.3)(.3,.3)
\rput(.9,.1){$\dots$}
\rput(1.3,.1){\tiny r}
\rput(.4,.1){\tiny 1}
\rput(.6,.1){\tiny 2}
\end{pspicture}
+
\begin{pspicture}[.2](-.2,0)(1.4,.6)
\psline{-}(0.3,0.1)(.3,.5)
\psline{-}(0.5,0.1)(.5,.5)
\psline{-}(1.2,0.1)(1.2,.5)
\psline[linestyle=dashed,dash=3pt 2pt]{-*}(.1,.3)(.5,.3)
\rput(.9,.1){$\dots$}
\rput(1.3,.1){\tiny r}
\rput(.4,.1){\tiny 1}
\rput(.6,.1){\tiny 2}
\end{pspicture}
+\dots
+\begin{pspicture}[.2](-.1,0)(1.4,.6)
\psline{-}(0.3,0.1)(.3,.5)
\psline{-}(0.5,0.1)(.5,.5)
\psline{-}(1.2,0.1)(1.2,.5)
\psline[linestyle=dashed,dash=3pt 2pt]{-*}(.1,.3)(1.2,.3)
\rput(.9,.1){$\dots$}
\rput(1.3,.1){\tiny r}
\rput(.4,.1){\tiny 1}
\rput(.6,.1){\tiny 2}
\end{pspicture}
$$
The above local image is called an {\em r-duplicated vertex.\/}
Applying Lemma~\ref{lemcom} in the case where the intersection of $D$ and the involved diagrams is an element of $\CA(\tribr)$ yields the following lemma:

\begin{lemma}
\label{lemcomd}
 The elements of $\CA(\tribr)$ commute with r-duplicated vertices.
\end{lemma}

\section{Good functors from q-tangles to $\CA$}
\label{secfunc}
\setcounter{equation}{0}

\begin{definition}
A {\em non-associative word\/} or {\em n.a. word} $w$ in the letter $\cdot$ is an element of the
free non-associative monoid generated by $\cdot$. The {\em length\/}
of such a $w$ is the number of letters of $w$. Equivalently, we can define
a {\em non-associative word \/} by saying that each such word has an integral {\em length\/} $\ell(w) \in \NN$, the only word of length 0 is the {\em empty
word\/}, the only word of length 1 is $\cdot$, the product $w^{\prime}w^{\prime \prime}$ of two
n.a. words $w^{\prime}$ and $w^{\prime \prime}$ is a n.a. word 
of length $\left(\ell(w^{\prime}) + \ell(w^{\prime \prime})\right)$, and every word $w$ of length
$\ell(w) \geq 2$ can be decomposed in a unique way as the product $w^{\prime}w^{\prime \prime}$ of two
n.a. words $w^{\prime}$ and $w^{\prime \prime}$ of nonzero length.
\end{definition}

\begin{example}
The unique n.a. word of length 2 is $(\cdot \cdot)$. The two
n.a. words of length 3 are $((\cdot \cdot) \cdot)$ and  $(\cdot (\cdot \cdot))$.
There are five n.a. words of length 4.
\end{example}

In the {\em ambient space\/} $\RR^3=\{(x,y,z)\}$, the {\em horizontal plane\/} is the plane $(z=0)$, whereas the {\em blackboard plane\/} is the plane $(y=0)$. The $z$-coordinate
of a point $(x,y,z) \in \RR^3$ is called its {\em vertical projection\/}.

\begin{definition}

A {\em q-tangle\/} is a triple $(T(M);b,t)$ where
$b$ and $t$ are two non-associative words and $T$ is a $C^{\infty}$ embedding of a compact one-manifold $M$ into a horizontal slice $\RR^2 \times [\beta,\tau]$ of $\RR^3$ such that:
$$T(M) \cap (\RR^2 \times \{\beta,\tau\}) = T(\partial M) \subset \RR \times \{0\} \times \{\beta, \tau\}$$
and the set of letters of $b$ and $t$ are in natural one-to-one correspondences
induced by the order of $\RR$
with $T^{-1}(\RR \times (0,\beta))$ and $T^{-1}(\RR \times (0,\tau))$,
respectively.
The only horizontal tangent vectors of $M$ occur for interior points of $M$ and are parallel to the blackboard plane. 
$T$ is considered up to the isotopies which satisfy these hypotheses at any time and up to a {\em rescaling of the height parameter\/} that is a composition by $1_{\RR^2}
\times h$ where $h$ is an increasing diffeomorphism from $[\beta,\tau]$ to another interval of $\RR$. The letters $b$ and $t$ stand for the 
{\em bottom word\/} and the {\em top word\/}, respectively.
Here, q-tangles will be simply called {\em tangles.}
\end{definition}

\begin{remark}
\label{rklim}
The involved non-associative words represent limit configurations of distinct points on the real line, that are corners of a suitable compactification
of the quotient of $\{(x_1, x_2, \cdots, x_p) \subset \RR^p / 
x_1 < x_2 < \cdots < x_p\}$ by the translations which identify
$(x_1, x_2, \cdots, x_p)$ to $(x_1 + T, x_2 + T, \cdots, x_p + T)$ for all $T \in \RR$ and by the positive homotheties which identify
$(x_1, x_2, \cdots, x_p)$ to $(\lambda x_1, \lambda x_2, \cdots, \lambda x_p)$
for all $\lambda > 0$. This compactification and the description of its strata may be found in \cite[Subsection 10.1]{p3}.
The points inside some matching parentheses are infinitely closer to
each other that they are to points outside the parentheses.
Here, it is enough to keep this interpretation in mind at an intuitive level.
\end{remark}
\begin{remark}
Each component $T(C)$ of a tangle (that is each image under $T$ of a connected component $C$ of
the underlying one-manifold $M$), has a well-defined {\em writhe\/} $w(T(C))$,
which is the sum of the signs of the self-crossings of $T(C)$ in the (regular
for a well-chosen $T$) projection onto the blackboard plane:
$$w(T(C))= (\mbox{number of}\; \fcpoc) - (\mbox{number of}\; \fcmoc).$$
To compute the writhe, we choose an arbitrary orientation of $C$, reversing it
does not affect the result.
\end{remark}

\begin{examples}
Tangles are unambiguously defined by the data of a regular projection of the involved embedding onto the
blackboard plane, together with the bottom and top words.
Since there is only one n.~a. word of length 0, 1 or 2, these words do not need to be specified. Here is an example of a q-tangle:
$$\smin =\left( \smin ;\emptyset,(\cdot \cdot ) \right) \neq \sminc$$
With Remark~\ref{rklim} in mind, the bottom and top words are sometimes shown in pictures by the relative positions of the bottom (or top) points. For example,
$$\ass=\left(\trit ;((\cdot \cdot ) \cdot),( \cdot( \cdot \cdot )\right)$$
$$\dcap=\left(\dcap ;((\cdot \cdot )(\cdot \cdot )),\emptyset \right)$$
\end{examples}

\noindent We can define the following operations on q-tangles: 

\begin{enumerate}
\item
The {\em product\/} of two q-tangles $T_1=(T_1(M_1);b_1,t_1)$ and $T_2=(T_2(M_2);b_2,t_2)$
is defined as soon as $t_1=b_2$ as the tangle $T=T_1T_2=(T(M);b_1,t_2)$ obtained by {\em stacking $T_2$ above $T_1$\/}:
These tangles are represented by embeddings $T$, $T_1$ and $T_2$
such that there exists a regular value $\gamma$ of the vertical
projection of the embedding $T:M \longrightarrow \RR^2 \times [\beta,\tau]$ such that $M_1=T^{-1}(\RR^2 \times [\beta,\gamma])$, $M_2=T^{-1}(\RR^2 \times [\gamma,\tau])$, and $T_1$ and $T_2$ are the restrictions of $T$ to $M_1$ and $M_2$, respectively. For example,
$$ \left(\smin\right) \left(\fcm\right)=\sminc$$
\item The {\em tensor product\/}  of two tangles $T_1=(T_1(M_1);b_1,t_1)$ and $T_2=(T_2(M_2);b_2,t_2)$
is defined as the tangle $T=T_1\otimes T_2=(T_1 \otimes T_2(M_1 \coprod M_2); b_1b_2, t_1t_2)$ by {\em putting $T_2$ on the right-hand side of $T_1$}.
In order to construct the embedding $T$, choose representatives of $T_1$ and $T_2$ which embed $M_1$ and $M_2$ into $[0,1] \times \RR \times [0,1]$ and 
$[2,3] \times \RR \times [0,1]$, respectively, and define $T$ as their disjoint union. For example,
$$\fcp \otimes \trid=\left(\fcp \trid;((\cdot\cdot)(\cdot\cdot)),((\cdot\cdot)(\cdot\cdot))\right)$$
\item The {\em duplication\/} of a component $C$ of a tangle $T$ consists in replacing $T(C)$ by two closed parallel copies
$T(C_1)$ and $T(C_2)$ of $T(C)$ so that, up to homotopy with fixed boundary, the section in the normal unit bundle of $T(C_1)$ induced by $T(C_2)$ coincides
with one of the two sections given, thanks to the condition on horizontal tangencies, by "the intersection with the blackboard plane". Every letter of the top and bottom words corresponding
to a (possible) boundary point of $C$ is replaced by the two-letter word
$(\cdot \cdot)$. The resulting tangle is denoted by $(2 \times C)(T)$. For example, duplicating the unique component of $\smin$ yields
$\dcup$. As another example, duplicating a knot $T(C)$ amounts to replace it by two parallel copies of $T(C)$ whose linking number is $w(T(C))$.
\item The {\em deletion\/} of a component $C$ of $M$ consists in forgetting
about $C$, removing the possible letters of the top and bottom words corresponding to the boundary points
of $C$ and removing the unneeded parentheses. The tangle obtained from $T=T(M;b,t)$ by deleting $C$ will be denoted by $T\setminus C$.
\item
The {\em orthogonal symmetry $s_h$ with respect to the horizontal plane \/} acts on tangles by a composition of the involved embedding by $s_h$ and by exchanging the top and bottom words.
The {\em orthogonal symmetry $s_v$ with respect to the blackboard plane  \/} acts on tangles by a composition of the involved embedding by $s_v$.
The
{\em  180-degree rotation $r_v$ around the vertical axis \/} acts on tangles by a composition of the involved embedding by $r_v$,
and by flipping both the top word and the bottom word.
\end{enumerate}

\begin{definition}
A {\em functor} from the category of q-tangles to $\CA$ is a map $Z$ 
which associates an element  $Z(T(M);b,t) \in  \CA(M)$ to any q-tangle $T=(T(M);b,t)$ so that $Z$ is compatible with the products, that is such that, if 
$T_1$ and $T_2$ are as in the definition of the products
on tangles, we have 
$$Z(T_1T_2)=Z(T_1)Z(T_2)$$
Graphically, this reads
$$ Z\left(
\begin{pspicture}[.4](0,0)(.7,1.4)
\psframe(.05,.7)(.65,1.3)
\psframe(.05,.1)(.65,.7) 
\rput(.35,1){$T_2$}
\rput(.35,.4){$T_1$}
\end{pspicture}
\right) =
\begin{pspicture}[.4](0,0)(1.4,1.4)
\psframe(.05,.7)(1.35,1.3)
\psframe(.05,.1)(1.35,.7) 
\rput(.7,1){$Z(T_2)$}
\rput(.7,.4){$Z(T_1)$}
\end{pspicture}$$

Such a functor is said to be {\em monoidal} if it respects tensor products,
that is if 
$$Z(T_1 \otimes T_2)=Z(T_1) \otimes Z(T_2)$$
Graphically, this reads
$$ Z\left(
\begin{pspicture}[.2](0,0)(1.3,.7)
\psframe(.65,0)(1.25,.6)
\psframe(.05,0)(.65,.6) 
\rput(.95,.3){$T_2$}
\rput(.35,.3){$T_1$}
\end{pspicture}
\right) =
\begin{pspicture}[.2](0,0)(2.7,.7)
\psframe(.05,0)(1.35,.6)
\psframe(1.35,0)(2.65,.6) 
\rput(2,.3){$Z(T_2)$}
\rput(.7,.3){$Z(T_1)$}
\end{pspicture}$$

Such a functor is said to {\em be compatible with the duplication of a regular
component\/} if for any component $C$ of a tangle $T$ which can be represented without horizontal tangent vector, we have 
$$Z((2 \times C)(T))=(2 \times C)_{\ast}(Z(T))$$

\end{definition}

Let $C$ be a component of a one-manifold $M$. We have a natural
map $\CO_C$ from $\CA(M)$ to $\CA(M \setminus C)$ such that:

$\CO_C([\Gamma]) =[\Gamma]$ if $\Gamma$ is a diagram without any leg on $C$ 
and $\CO_C([\Gamma]) =0$ for the other diagrams.

A functor $Z$ is said to {\em respect deletion\/} (or to be compatible with deletion) if for any component $C$ of a tangle $T$, we have 
$$Z(T\setminus C)=\CO_C(Z(T))$$

A functor $Z$ is said to be {\em invariant under the 180 degree rotation $r_v$ around
the vertical axis\/} if, for any q-tangle $T$, $Z \circ r_v(T)=Z(T)$

\begin{definition}
An element of $\CA_n$ is said to be {\em real (or rational)\/} when it can be written as a combination of diagrams with real (or rational) coefficients.
Equivalently, define $\CA_n^{\bf Q}(M\cup X)$ as the rational vector space 
generated by the diagrams with support $M \cup X$ quotiented out by the
relations $AS$, $STU$ and $IHX$. Then  
$\CA_n(M\cup X) = \CA_n^{\bf Q}(M\cup X) \otimes_{\QQ} \CC $. Set
$\CA_n^{\bf R}(M\cup X) = \CA_n^{\bf Q}(M\cup X) \otimes_{\QQ} \RR $,
then the rational elements of $\CA_n$ are the elements of $\CA_n^{\bf Q}$ whereas
the real elements of $\CA_n$ are the elements of $\CA_n^{\bf R}$.
An element of $\CA$ is said to be {\em rational\/} if it lies in $\prod_{n \in \NN}\CA^{\bf Q}_n \subset \prod_{n \in \NN}\CA_n$.
\end{definition}

It follows from \cite{p3} that the Poirier limit integral is a good monoidal functor which, in addition, is real for all tangles and rational for framed links.
It follows from \cite{lm} and the remark below that the Kontsevich integral is also a good monoidal functor with $a^{Z^K}=\frac12\tata$. See also \cite{les}.

\begin{remark}
\label{rkor}
Our definition of diagrams on manifolds $M$ differs from the one used in \cite{lm}
or \cite{bn} where the univalent vertices are not oriented but the manifolds are.
Nevertheless, note that the datum of the orientation of a univalent vertex is equivalent to the datum of a local orientation of $M$ near the univalent vertex. Namely, we identify  
\begin{pspicture}[.2](0,0)(.8,.4)
\psline[linestyle=dashed,dash=3pt 2pt]{-*}(.35,.4)(.35,.1)
\psline{-}(.1,.1)(.7,.1)
\end{pspicture} (without any specified orientation on the solid line) to  \begin{pspicture}[.2](0,0)(.8,.4)
\psline[linestyle=dashed,dash=3pt 2pt]{-*}(.35,.4)(.35,.1)
\psline{->}(.1,.1)(.7,.1)
\end{pspicture}. Then the antisymmetry relation (AS) for univalent vertices allows us to represent all the diagrams by linear combinations of diagrams
where the orientation of univalent vertices match the global orientation of manifolds. These diagrams are the only ones considered in \cite{lm} or \cite{bn} (where there is no AS relation for univalent vertices).
\end{remark}

Now, we are in a position to understand the hypotheses of the theorem.
Let us begin with some known remarks in order to understand its proof.

\section{Classical properties of good monoidal functors}
\setcounter{equation}{0}

\begin{definition}
The symmetry $\sigma_f$ of $\CA(\triu)$ with respect to the flip of the interval is the symmetry induced by the flip of $\triu$: the injection of the legs is composed by the flip which also carries (and thus changes) the local orientations of the legs.
It is unknown whether there exist elements of $\CA(\triu)$ that are not symmetric under this symmetry.

\end{definition}
\begin{proposition}
\label{propgood}
Let $Z$ be a good monoidal functor. \\Let $\Phi_Z=Z(\ass) \in \CA(\trit)$. \\
Let 
$\nu(Z) \in \CA(\triu)$ be the element such that $\nu(Z)_0=1$, and the map from $\CA(\ass)$ to $\CA(\bigsnake)$ induced by the inclusion shown by the pictures maps $\Phi_Z$ to $\nu(Z)^{-2}$.\\
Then $Z$ satisfies the properties below.
Furthermore, they remain true if $Z$ is replaced by its truncature at order $k$ that is the composition of $Z$ by the quotient map from $\CA$ to $\frac{\CA}{\prod_{r\geq k+1} \CA_r}$.
\begin{enumerate}
\item $Z(\triu)=1$. 
\item $Z(\fcp) \in \CA(\trid)$ is symmetric
with respect to the exchange of the two strands, and, in $\CA(\trid)$, 
$$Z(\fcm)=Z(\fcp)^{-1}.$$
\item $Z(\smin)$ and $Z(\smax)$ are symmetric elements of $\CA(\triu)$ with respect to the flip of $\triu$ and 
\begin{equation}\label{eqnu}Z(\smin)Z(\smax)=\nu(Z)^{2}=Z(\scirc).\end{equation}
\item $a^Z \in \CA(\triu)$ is symmetric with respect to the flip of $\triu$.
\item For any q-tangle $T=(T(M);b,t)$, $Z(T)$ only depends on $b$, $t$, the (unframed -the horizontal tangent vectors are allowed to turn during the isotopy-) isotopy class of $T$ (with fixed boundary order) and on the writhe of each component.
\item The elements $Z(\dcupn)$ and $Z(\dcapn)$ of $\CA(\tridn)$ are symmetric with respect to the simultaneous flip of the two components, and we have, in $\CA(\tridn)$,
\begin{equation}
\label{eqdup}
Z(\dcupn) Z(\dcapn) = \left( 2 \times \triu\right)_{\ast}(Z(\smin)Z(\smax)) = \left( 2 \times \triu\right)_{\ast}(\nu(Z)^{2}).\end{equation}
\item 
\begin{equation}\label{eqtwist}Z(\fcp)= \exp(\frac12(2 \times \triu)_{\ast} a^Z)\left(\exp(-\frac12a^Z)\otimes \exp(-\frac12a^Z)\right) \;\;\; \mbox{in} \; \CA(\trid).\end{equation}
\item The elements $Z(\dcupn)$ and $Z(\dcapn)$ of $\CA(\tridn)$ are symmetric with respect to the exchange of the two components of $\tridn$.
\item $Z$ is compatible with the duplication of any component which has as
many minima and maxima.
\item On tangles, $Z$ is determined by $\Phi_Z$, $a^Z$ and $Z(\smax)$.
\item On framed links, $Z$ is determined by 
$\Phi_Z$ and $a^Z$.
\item $\nu(Z)$ is real and is zero in odd degrees.
\item For any integer $i$, $a^Z_{2i}=0$ and $a^Z_{2i+1}$ is real.
\end{enumerate}
\end{proposition}
\bp We proceed with respect to the order of the statement.
\begin{enumerate}
\item Since $Z(\triu)^2=Z(\triu)$ and $Z_0(\triu)=1$, 
$$Z(\triu)=1.$$
\item It follows by duplication that $Z(\trid)=1 \in \CA(\trid)$.
Now, by functoriality, and because of the conditions on the degree zero parts,
$$Z(\fcm)=Z(\fcp)^{-1}$$ in $\CA(\trid)$ where $Z(\fcp)$ must be symmetric
w.r.t. the exchange of the two strands because of the invariance under the 
180-degree rotation $r_v$ around the vertical axis.
\item 
The invariance of $Z$ under $r_v$
implies that $Z(\smax)$ and $Z(\smin)$ are symmetric w.r.t. the flip of the interval. Now, we prove $Z(\smin)Z(\smax)=\nu(Z)^{2}$ by computing $Z(\snake)$.
The properties of $Z$ ensure that: 
$$Z(\snake)=Z(\smin)Z(\smax)\nu(Z)^{-2}$$
The above symmetry property of $Z(\smax)$ and $Z(\smin)$ allows us to
forget to specify an orientation for $\snake$ to define the $\CA(I)$-module structure involved.
Now, we get Equation \ref{eqnu} because $\snake$ is isotopic as a q-tangle to $\triu$.
\item 
$a^Z \in \CA(\triu)$ is symmetric with respect to the flip of $\triu$, because both
$Z(\smax)$ and $Z(\smaxn)$ are invariant under $r_v$.
\item 
If two framed tangles are related by an unframed isotopy, then they are related
by a sequence of framed isotopies and twists of the extrema. During such a sequence, $Z$ changes by a multiplication of $\exp((\delta w) a^Z)$ on each component, where $\delta w$ represents the variation of the writhe on the given component.
\item 
The invariance of $Z$ under $r_v$ implies the symmetry of $Z(\dcupn)$ and $Z(\dcapn)$  w.r.t. the simultaneous flip of the two intervals. Now, the computation of $Z(2 \times \snake)=Z(\trid)=1$
and the fact that duplicated vertices commute with the elements of $\CA(\trid)$
(Lemma~\ref{lemcomd}) implies Equation~\ref{eqdup} as follows:
$$1=Z(2 \times \snake)=\left(Z(\dcupn)\otimes 1_{\itrid}\right)
\left((2 \times \ass)_{\ast}Z(\ass) \right)
\left( 1_{\itrid} \otimes Z(\dcapn) \right)$$
$$= \left( (2 \times \triu)_{\ast} \left( \nu(Z)^{-2} \right) \right)
\left( Z(\dcupn)Z(\dcapn) \right)$$
\item
Now, we determine $Z(\fcp)$ as a function of $a^Z$ as in \cite{p3}
using the belt trick, by computing $Z$ for the tangle obtained from $\hbbelt$
by duplication.
Indeed, $$Z(\hbbelt)=\exp(2a^Z)Z(\snake)=\exp(2a^Z)$$
whereas 
Equation~\ref{eqdup}, Lemma~\ref{lemcomd} and the duplication property for regular strands yield:
$$Z(2\times \hbbelt)=(2 \times \triu)_{\ast}Z(\hbbelt)$$
Thus, $$Z(2\times \hbbelt)=(2 \times \triu)_{\ast}(\exp(2a^Z))=\exp((2 \times \triu)_{\ast} (2a^Z))$$
On the other hand, the q-tangle $2 \times \hbbelt$ is equal to the product
$\left(\fcp \right)^4 \left(\hbbelt \otimes \hbbelt \right)$. (We could also use the above stronger invariance property (5) of $Z$.) Thus, 
$$Z(2\times \hbbelt)=Z(\fcp)^4\left(\exp(2a^Z)\otimes \exp(2a^Z)\right)$$
and 
$$Z(\fcp)^4=\exp((2 \times \triu)_{\ast} (2a^Z))\left(\exp(-2a^Z)\otimes \exp(-2a^Z)\right)$$
This gives the announced formula for $Z(\fcp)$ because the terms in this formula commute in $\CA(\trid)$. 
\item
Now, we can just use the same proof in the opposite direction to deduce that
$Z(2\times \hbelt)=(2 \times \triu)_{\ast}Z(\hbelt)$ from the expression of  $Z(\fcp)$, and we may deduce from this fact that
$$Z(\dcupn) Z(\dcapnof) = \left( 2 \times \triu\right)_{\ast}(Z(\smin)Z(\smax))$$
Thus, using Equation~\ref{eqdup},
$$ Z(\dcapn) = Z(\dcapnof) = Z(\dcapno)$$
Similarly, $Z(\dcupn)$ is symmetric with respect to the exchange of the two intervals, too.
\item
Now, this symmetry, Equation~\ref{eqdup} and Lemma~\ref{lemcomd} imply that $Z$ is compatible with the duplication of the components which have as many minima and maxima.
\item
Every q-tangle can be expressed as a product of q-tangles which are
either tensor products of elements of the form $\fcm$, $\fcp$, $\smax$, $\smin$  and $\triu$, or q-tangles of the form
$(\trib;b,t)$ which in turn can be expressed as products of tensor products
of elements that can be obtained by duplicating $\ass$, $r_v\left(\ass\right)$ and $\triu$.
Thus, the relations that are already shown prove that $Z$ is uniquely determined by
$a^Z$, $\Phi_Z$ and $Z(\smax)$.
\item
Since each component of a link has as many maxima and minima, knowing $Z(\smax)$ is
unneeded to determine $Z$ on links.
\item
Now, note that $Z(\scirc)=Z(\smin)Z(\smax)$ is unchanged by the symmetries $\sigma_v$ and
$\sigma_h$. Thus, $\nu(Z)^{2}$ is zero in odd degrees and is real in even degrees, and 
the same must be true for $\nu(Z)$.
\item
Since \huit is obtained from \huito by the symmetries $s_v$ and $s_h$, we have: 
$$\exp(-a^Z)=\sigma_h(\exp(a^Z))=\sigma_v(\exp(a^Z))$$
This proves the last statement of the proposition.
\end{enumerate}
\eop

\begin{definition}
An element of $\CA(M\cup X)$ is {\em of filtration at least d\/} if it is a combination of diagrams of degree at least d.
\end{definition}
	
In \cite{lm}, Le and Murakami proved that, if $Z$ is a good 
monoidal functor such that  $a^Z=a^{Z^K}=\frac12\tata$, then $Z$ coincides
with the Kontsevich integral on framed links. More precisely, and that is what will be needed for our proof, they proved the following result:

\begin{theorem}[Le-Murakami \cite{lm}]
\label{thlm}
If $Z$ is a good monoidal functor such that  $\left(a^Z-a^{Z^K}\right)$ is of filtration at least $(k+1)$, then for any framed link $L$, $\left(Z(L)-Z^K(L)\right)$ is of filtration at least $(k+1)$.
\end{theorem}

Since this result is not stated in these words in \cite{lm}, we will review
its proof in Section~\ref{seclm}.

\section{More algebra on diagrams}
\label{secmad}
\setcounter{equation}{0}

This section is mostly devoted to stating and (re-)proving Lemma~\ref{lemeig} that will be used in the proof of Theorem~\ref{mainth}. The arguments there
have already been used by D.~Bar-Natan, T.~Le and D.~Thurston to compute the Kontsevich Integral of the unknot. I have learned them from Dylan Thurston.

Let $x \in X$. 
We now describe the so-called Poincar\'e-Birkhoff-Witt isomorphism $\chi_x$
from $\CA(M \cup X)$ to $\CA(M \cup I \cup  (X \setminus \{x\}))$. 

Let $\Gamma$ be a diagram with support $M \cup X$ as in Definition \ref{defdia}. Let $U_x$ denote
the preimage of $\{x\}$ under $f$ and let $k$ be its cardinality.
Let $\Sigma$ be the set of isotopy classes of injections $i$ from $U_x$ into the interior of $I$. $\Sigma$ is in one-to-one correspondence with the set of the $k!$ total orders on $U_x$. Every element $\sigma$ of $\Sigma$ naturally defines an oriented diagram $\Gamma_{\sigma}$ with support $(M \cup I \cup  (X \setminus \{x\}))$ where the local orientations 
of the elements of $U_x$ are induced by the usual orientation of $I$ as in Remark~\ref{rkor}.
Let $\chi_x$ be the morphism of topological vector spaces which maps the class $[\Gamma]$ of a diagram $\Gamma$ as above to:
$$\chi_x([\Gamma])= \frac1{k!} \sum_{\sigma \in \Sigma}[\Gamma_{\sigma}]$$
It is easy to see that $\chi_x$ is well-defined and the proof of \cite[Theorem 8]{bn} applies to prove that $\chi_x$ is an isomorphism. (The surjectivity of 
$\chi_x$ is reproved along the same lines in the proof of Lemma~\ref{lempbwo} below.)

The spaces $\CA(\{1, \cdots,r\})$ and $\CA(\{1\})$ are simply denoted by
$\CB(r)$ and $\CB$, respectively.
The (commutative) product $\chi_1 \circ \chi_2 \circ \cdots \circ \chi_r$ from $\CB(r)$
to $\CA(\{1, \cdots,r\} \times I)$ is denoted by $\chi(r)$.

We have the following easy sublemma.
\begin{sublemma}
\label{subu}
The following diagram is commutative:
$$\diagram{\CB&\hfl{\chi}&\CA(I) \cr\vfl{(r \times 1)_{\ast}}&&\vfl{(r \times I)_{\ast}}
\cr \CB(r)&\hfl{\chi(r)}&\CA(\{1, \cdots,r\} \times I)}$$
\end{sublemma}

Let $\CB_{n,k}$ be the subspace of $\CB_n$ generated by the degree n diagrams with support $\{1\}$
with exactly $k$ legs, and let $\CB_{\ast,k}=\prod_{n \in {\bf N}}\CB_{n,k}$.
With this notation, an element $\beta$ of $\CA(I)$ has two legs if and only if 
$\chi^{-1}(\beta) \in \CB_{\ast,2}$. 

Let $j$ denote the linear continuous map from $\CB(r)$ to $\CB$  which maps
a diagram $\Gamma$ with support $\{1,\cdots, r\}$ as in Definition~\ref{defdia} to the diagram
with support $\{1\}$ obtained by replacing its function $f$ by the constant
function. 

Note the following sublemma.
\begin{sublemma}
\label{subuu}
Let $k$ and $r$ be two integers, $r \geq 2$, then $\CB_{\ast,k}$ is the eigenspace of the endomorphism
$j \circ (r \times 1)_{\ast}$ associated to the eigenvalue $r^k$.
\end{sublemma}

Consider the map $\iota$ from $\{1, \cdots,r\} \times I$ to $I$ which maps $(k,t)$ to $\frac{k-1+t}r$, it naturally induces a map $\iota_{\ast}$ from
$\CA(\{1, \cdots,r\} \times I)$ to $\CA(I)$.

The proof of the following sublemma is again left as an exercise for the reader.
\begin{sublemma}
\label{subuuu}
The following diagram is commutative:
$$\diagram{\CB&\hfl{(r \times 1)_{\ast}}& \CB(r) & \hfl{j} & \CB \cr
\vfl{(r \times 1)_{\ast}}&&&&\vfl{\chi} \cr
\CB(r)&\hfl{\chi(r)}&\CA(\{1, \cdots,r\} \times I)&\hfl{\iota_{\ast}}&\CA(I)}$$
\end{sublemma}

The following known\footnote{In \cite{bn}, $\iota_{\ast} \circ (r \times I)_{\ast}$
is called the $r^{th}$ Adams operation.} lemma is an easy consequence of the three previous sublemmas.

\begin{lemma}
\label{lemeig}
The eigenspace of the endomorphism 
$$ \CA(I) \; \hfl{(2 \times I)_{\ast}} \; \CA(\{1,2\} \times I) \; \hfl{\iota_{\ast}} \CA(I)$$ 
associated to the eigenvalue $4$ is exactly the subspace $\chi(\CB_{\ast,2})$
of $\CA(I)$ made of the two-leg elements of $\CA(I)$.
\end{lemma}

\section{The isomorphisms $\Psi(\beta)$}
\setcounter{equation}{0}

Let $\beta$ be a two-leg element of $\CA(I)$. 
Let $\Gamma$ be a diagram with support $M \cup X$ as in Definition~\ref{defdia}. We define $\Psi(\beta)(\Gamma)$ to be the element of $\CA(M \cup X)$ obtained by inserting $\beta^s$ $d$ times on each degree $d$
component of $\Gamma$ (where a {\em component\/} of $\Gamma$ is a connected component of the dashed graph). See Definition~\ref{defbetas}.

\begin{lemma}
\label{lemins}
$\Psi(\beta)(\Gamma)$ does not depend on the choice of the insertion loci.
\end{lemma}
 
\bp
It is enough to prove that moving $\beta^s$ from an edge of $\Gamma$ to another one
does not change the resulting element of $\CA(M \cup X)$, when the two
edges share some vertex $v$. Since this move amounts to slide $v$ through 
$\beta^s$, it suffices to prove that sliding a vertex from some leg of a two
leg-diagram to the other one does not change the diagram modulo $AS$, $IHX$ and $STU$, this
is a direct consequence of Lemma~\ref{lemcom} when the piece of diagram
inside $D$ is $\beta^s$.
\eop

It is now easy to check that $\Psi(\beta)$ is compatible with the relations IHX, STU and AS. This allows us to define continuous vector space endomorphisms  $\Psi(\beta)$ of 
the $\CA(M \cup X)$ such that, for any diagram $\Gamma$:
$$\Psi(\beta([\Gamma]))= \Psi(\beta)(\Gamma)$$

$\Psi(\beta)$ satisfies the following properties:
\begin{lemma}
\label{lemPhi}
\begin{enumerate}
\item $\Psi(\beta)$ is compatible with the products of Section~\ref{secdia}.
($\Psi(\beta(xy))=\Psi(\beta(x))\Psi(\beta(y))$.)
\item  $\Psi(\beta)$ commutes with the duplication maps of Section~\ref{secdia}.
\item $\Psi(\beta)$ commutes with the Poincar\'e-Birkhoff-Witt isomorphisms
$\chi_x$.
\item If $\beta_1 \neq 0$, $\Psi(\beta)$ is an isomorphism of topological vector spaces such that $\Psi(\beta)$ and $\Psi(\beta)^{-1}$ map elements of filtration
at least $d$ to elements of filtration at least $d$.
\item If $\beta_1 \neq 0$, if $\beta$ is real and null in even degrees, and if $Z$ is a good monoidal functor from the category of q-tangles to $\CA$, then  $\Psi(\beta) \circ Z$ and $\Psi(\beta)^{-1} \circ Z$ are good monoidal functors, too.
\end{enumerate}
\end{lemma}
\bp
The first three properties are obvious. For the fourth one, first note that 
$\beta_1 = b_1 \tata$ for some non zero number $b_1$. Thus, for $x =\sum_{i=d}^{\infty}x_i$, $\Psi(\beta)(x) - (b_1^d x_d)$ is of filtration at least $d+1$. This shows that $\Psi(\beta)$ is injective and allows us to construct 
a preimage for any element by induction on the degree, proving that $\Psi$ is onto.
The fifth property is a consequence of the other ones.
\eop

\section{Proof of the theorem}
\setcounter{equation}{0}

Set $a^Z=a$. $a_1\neq 0$, and, by Proposition~\ref{propgood}, we know that $a_{2n}$ is zero for any integer $n$.
We proceed by induction on $N$, and we work with the induction hypotheses:

\medskip

\noindent
$(\star(N))$ $a_k$ is a two-leg element for any $k \leq 2N$.
Let $\Psi_N$ denote the isomorphism $\Psi(2 \sum_{i=1}^{2N}a_i)$.

\medskip

\noindent
$(\star \star(N))$ For any framed link $L$, $\left(\Psi_N \circ Z^K - Z\right)(L)$ is of filtration at least $2N+1$.

\medskip

\noindent
Of course, $(\star(1))$ is true.
We are going to prove the two lemmas:

\begin{lemma}
\label{lemin1}
For any integer $N \geq 1$,  $(\star(N))$ implies $(\star \star(N))$.
\end{lemma}

\begin{lemma}
\label{lemin2}
For any integer $N \geq 1$,  $(\star(N))$ and $(\star \star(N))$ imply $(\star(N+1))$.
\end{lemma}

Assume that these lemmas are proved, then $a$ is a two-leg element of $\CA(I)$, and, for any framed link $L$, for any $n \in {\bf N}\setminus \{0\}$, the degree n parts of $\Psi(2a)(Z^K(L))$ and $\Psi_n(Z^K(L))$ coincide. Hence, the degree n part
of  $\Psi(2a)(Z^K(L))$ is equal to the degree n part of $Z(L)$
and the theorem is proved except for the two lemmas whose proofs follow.

\medskip

\noindent{\sc Proof of Lemma \ref{lemin2}:}
Let $Z$ be a good monoidal functor.
Set $$A(Z)= Z\left(\huit\right) = \exp(a^Z)Z(\scirc)$$
and $$B(Z)= Z\left(\DT\right) = Z(\dcup) \left(Z(\fcp) \otimes Z(\trid) \right) Z(\dcap)$$
Since $\huit$ and $\DT$ are both isotopic to the trivial unframed knot, and since both  have writhe 1, 
$A(Z)=B(Z)$. This reads $A_{2N+1}(Z)=B_{2N+1}(Z)$ in degree $(2N+1)$.

In particular, for our $Z$, we have the equality:
\begin{equation}
\label{eqab}
A_{2N+1}(Z)-A_{2N+1}(\Psi_N \circ Z^K)=B_{2N+1}(Z)-B_{2N+1}(\Psi_N \circ Z^K).
\end{equation}
Since $Z(\scirc) =\nu(Z)^2$ is null in odd degrees according to
Proposition~\ref{propgood}, and because $(\star \star(N))$ is assumed to be true, $\left(Z(\scirc)-\Psi_N \circ Z^K(\scirc)\right)$
is of filtration at least $(2N+2)$. Furthermore, 
$$a^{\Psi_N \circ Z^K}=\Psi_N(a^{Z^K})=\sum_{i=1}^{2N}a_{i}.$$
Therefore, the
left-hand side of \ref{eqab} is equal to $a_{2N+1}$.

Similarly, using Equation~\ref{eqdup}, the symmetries (6) and (8) of Proposition~\ref{propgood} and $(\star \star(N))$,
we can see that
$$Z(\dcupn)Z(\dcapnof)-\Psi_N \circ Z^K(\dcupn)\Psi_N \circ Z^K(\dcapnof)$$
is of filtration at least $(2N+2)$.
Equation~\ref{eqtwist} shows us that $\left(Z(\fcp)-(\Psi_N \circ Z^K)(\fcp)\right)$ is of filtration at least $(2N+1)$. Thus,
$$B_{2N+1}(Z)-B_{2N+1}(\Psi_N \circ Z^K)=\iota_{\ast}\left(Z_{2N+1}(\fcp)-(\Psi_N \circ Z^K)_{2N+1}(\fcp)\right)$$
where $\iota_{\ast}$ is defined right before Sublemma~\ref{subuuu}.
Now, using Equation~\ref{eqtwist} both for $Z$ and for $\Psi_N \circ Z^K$, we get
$$B_{2N+1}(Z)-B_{2N+1}(\Psi_N \circ Z^K)=-a_{2N+1} + \frac12 \iota_{\ast} \circ (2 \times \triu)_{\ast} (a_{2N+1})$$
Hence,
$$ \iota_{\ast} \circ (2 \times \triu)_{\ast} (a_{2N+1})=4a_{2N+1}$$
and this proves that $\ast(N+1)$ is true, thanks to Lemma~\ref{lemeig}.
 \eop

\noindent{\sc Proof of Lemma \ref{lemin1}:}
$\Psi_N^{-1} \circ Z$ is a good monoidal functor such that 
$a^{\Psi_N^{-1} \circ Z}=\Psi_N^{-1}(a^Z)$ and  $\left(\Psi_N^{-1}(a^Z)-a^{Z^K}\right)=\Psi_N^{-1}\left(a^Z-\Psi_N(a^{Z^K})\right)$ is of filtration at least 2N+1. Thus, the Le and Murakami theorem (\ref{thlm}) applies and proves the lemma.
\eop

The proof of Theorem~\ref{mainth} is now complete except that
 Theorem~\ref{thlm} is stated neither in these words nor with our conventions in \cite{lm}. Therefore and in order to make the uniqueness statement more specific, we will
sketch the proof of Theorem~\ref{thlm} and refine its statement in the section below.

\section{More about the Le and Murakami theorem}
\label{seclm}
\setcounter{equation}{0}

In this section, we reduce the proof of Theorem~\ref{thlm} to the proof of \cite[Proposition 3, Section 8]{lm} that is proved in \cite{lm}, following \cite{lm} and using the same notation as in \cite{lm}.
Namely, $\CP_r$ denotes\footnote{In fact, $\CP_r$ is denoted by $\CP_r \otimes \CC$ in \cite{lm}. Warning: In this section, the subscripts do not represent degrees of diagram. } $\CA(\tribr)=\CA(\{1,2, \dots,r\} \times I)$. For $i \leq r$,
$\Delta_i:\CP_r \rightarrow \CP_{r+1}$ denotes the duplication of the $i^{th}$ strand $\left(2 \times \left(\{i\} \times I\right) \right)_{\ast}$, and 
$\varepsilon_i:\CP_r \rightarrow \CP_{r-1}$ denotes the deletion of the $i^{th}$ strand ($\CO_{\{i\} \times I}$). We say that an element of $\CP_2$ is {\em symmetric\/} if it is symmetric under the permutation of the two strands.

\begin{definition}
\label{deftwist}
A {\em twist\/} is a symmetric element $F$ of $\CP_2$ such that $\varepsilon_1(F)=1$.
Such a twist gives rise to an $\CF(w) \in \CP_{\ell(w)}$ for every n.~a. word $w$
which is uniquely defined by induction on $\ell(w)$ as follows:
$\CF(\emptyset)=1$, $\CF(\cdot)=1$, and
$$\CF(w^{\prime}w^{\prime \prime})=
\left(\left(\ell(w^{\prime}) \times I_1 \right)_{\ast}\circ \left(\ell(w^{\prime \prime}) \times I_2 \right)_{\ast}(F)  \right)
\left(\CF(w^{\prime}) \otimes \CF(w^{\prime \prime})\right)$$
where $I_1=\{1\} \times I$ and $I_2=\{2\} \times I$.
\end{definition}
Observe that $\CF((\cdot\cdot))=F$ and that the two factors in the above product commute thanks to Lemma~\ref{lemcomd}.
\begin{definition}
\label{defftwist}
Let $Z$ be a good monoidal functor and let $F \in \CP_2$ be a twist. Define
a functor $Z^F$ by the formula
$$Z^F(T(M);b,t)= \CF(b)^{-1}Z(T(M);b,t)\CF(t)$$
that also reads:
$$ Z^F(T(M);b,t) =
\begin{pspicture}[.4](0,0)(3,2)
\psframe(.05,1.3)(2.95,1.9)
\psframe(.05,.7)(2.95,1.3)
\psframe(.05,.1)(2.95,.7) 
\rput(1.5,1.6){$\CF(t)$}
\rput(1.5,1){$Z(T(M);b,t)$}
\rput(1.5,.4){$\CF(b)^{-1}$}
\end{pspicture}$$
\end{definition}

\begin{proposition}
If $Z$ is a good monoidal functor and if $F \in \CP_2$ is a twist, then
$Z^F$ is a good monoidal functor that coincides with $Z$ on framed links
and such that $a^{Z^F}=a^Z$.
\end{proposition}
\bp
It is obvious that $Z$ and $Z^F$ are functors which coincide on framed links and
that $Z^F_0=1$. The symmetry of $F$ makes clear that $F1_{\ismax} \in \CA(\triu)$ is symmetric with respect to the flip of the interval.
$a^{Z^F}=a^Z$ because $Z^F(\ssmaxn)=Z(\ssmaxn) \left( F1_{\ismax} \right)=\exp(a^Z)Z(\smax) \left( F1_{\ismax} \right)=\exp(a^Z)Z^F(\smax)$ in $\CA(\triu)$.
Now, let us prove that $Z^F$ is monoidal. 

 $$Z^F\left((T_1;b_1,t_1)\otimes (T_2;b_2,t_2)\right)
=$$
$$\left(\ell(b_1) \times I_1 \right)_{\ast}\circ \left(\ell(b_2) \times I_2 \right)_{\ast}(F^{-1})
\left( Z^F(T_1) \otimes Z^F (T_2) \right)
\left(\ell(t_1) \times I_1 \right)_{\ast}\circ \left(\ell(t_2) \times I_2 \right)_{\ast}(F)$$
Applying Lemma~\ref{lemcom} successively to all the vertices $v$ of any diagram $\Gamma$ with support $\trid$ shows that$$\left(\ell(b_1) \times I_1 \right)_{\ast}\circ \left(\ell(b_2) \times I_2 \right)_{\ast}([\Gamma])\left( Z^F(T_1) \otimes Z^F (T_2) \right)$$
$$=\left( Z^F(T_1) \otimes Z^F (T_2) \right)\left(\ell(t_1) \times I_1 \right)_{\ast}\circ \left(\ell(t_2) \times I_2 \right)_{\ast}([\Gamma])$$ 
because the sum of duplications of a vertex $v$ may slide through $Z^F(T_1)$ or $Z^F(T_2)$.
This proves that 
$$Z^F(T_1\otimes T_2)=Z^F(T_1) \otimes Z^F (T_2)$$
Now, the following properties of the $\CF(w)$ are easy to obtain by induction on the length. 
For any n.~a. word $w$ of length $\ell$:
\begin{enumerate}
\item deleting the $i^{th}$ letter of $w$ transforms $\CF(w)$ into $\varepsilon_i(\CF(w))$,
\item changing the $i^{th}$ letter of $w$ into the two-letter word $(\cdot \cdot)=(\cdot_i\cdot_{i+1})$ transforms $\CF(w)$ into $\Delta_i(\CF(w))\CF((\cdot_i\cdot_{i+1}))$, and
\item flipping $w$ into $r_v(w)$ transforms $\CF(w)$ into $r_v(\CF(w))$ where $r_v$ acts on $\CP_{\ell}$ by flipping the order of the vertical intervals: $\CF(r_v(w))=r_v(\CF(w))$
\end{enumerate}
Once these properties are proved, $Z^F$ is compatible with the duplication of a regular component, because the extra factors of the form $\CF((\cdot_i\cdot_{i+1}))^{\pm 1}$ may slide
along duplicated strands, and thus cancel each other with the help
of Lemma~\ref{lemcomd}. It is clear that $Z^F$ satisfies all the other properties of a good monoidal functor. 
\eop

\noindent Note that if $F$ and $F^{\prime}$ are two twists, then $FF^{\prime}$ is
a twist and $$Z^{FF^{\prime}}=\left(Z^F\right)^{F^{\prime}}$$

\noindent Recall that $$\Phi_Z=Z\left(\ass=\left(\trit ;((\cdot \cdot ) \cdot)
,( \cdot( \cdot \cdot ))\right)\right).$$
Because of Proposition~\ref{propgood} (11), 
the proof of Theorem~\ref{thlm}
easily reduces to the proof of the following statement. 

{\em
If $Z$ is a good monoidal functor such that  $\left( a^Z-a^{Z^K} \right) $ is of filtration at least $(k+1)$, then there exists a twist $F \in \CP_2$ such that 
$\left(\Phi_{Z^F}-\Phi_{Z^K}\right)$ is of filtration at least $(k+1)$. \/}

In turn, by induction and because of the previous remark on the composition of twists, the proof of this statement reduces to the proof
of the following lemma:

\begin{lemma}
\label{lemklm}
Let Z be a good monoidal functor such that $\left(a^{Z}-a^{Z^K}\right)$ is of filtration at least $k+1$ and such that
$\left(\Phi_Z -\Phi_{Z^K}\right)$ is of filtration at least $k$ ($k \geq 1$), then there exists a  symmetric degree $k$
element $f \in \CP_2$ such that $\varepsilon_1(f)=0$ and
$\left(\Phi_{Z^{(1+f)}} -\Phi_{Z^K}\right)$
is of filtration at least $(k+1)$.
\end{lemma}
\bp
Let $$d:\CP_n \rightarrow \CP_{n+1}$$ 
$$ f\mapsto d(f)=1 \otimes f - \Delta_1(f) + \Delta_2(f) - \dots +(-1)^n \Delta_n(f) + (-1)^{n+1} f \otimes 1$$
Now, let $f$ be a symmetric  degree $k$ element of $\CP_2 \cap \varepsilon_1^{-1}(0)$. By Definitions~\ref{deftwist} and~\ref{defftwist}, the part of $\left( \Phi_{Z^{(1+f)}} - \Phi_{Z} \right)$
of degree less than $(k+1)$ is equal to $d(f)$.
Thus, if $\psi$ denotes the degree $k$ part of $\left(\Phi_{Z^K}-\Phi_Z\right)$, it is
enough to find a symmetric, degree $k$ element of $\CP_2$ such that 
$\varepsilon_1(f)=0$ and 
$$d(f)=\psi$$
Now, what do we know about $\psi$? \\
First $\psi \in \CP_3$. When $\sigma$ is a permutation of $\{1,2,3\}$,
$\psi^{\sigma(1) \sigma(2) \sigma(3)}$ denotes the element of $\CP_3$ obtained from $\psi$ by sending the $i^{th}$ strand to the $\sigma(i)^{th}$ one. Then $\psi$ satisfies C1, C2, C3 and C4 below.
\begin{description}
\item[C1] $$d(\psi)=0$$
Indeed, we have the following {\em pentagon equation} in $\CP_4$
$$\Delta_1(\Phi_Z)\Delta_3(\Phi_Z)=\left(\Phi_Z \otimes 1_{\itriu}\right)\Delta_2(\Phi_Z)\left(1_{\itriu} \otimes \Phi_Z\right)$$
which is satisfied by both $\Phi_{Z^K}$ and $\Phi_Z$ and which expresses the fact that
$$Z\left(\triq;(((\cdot \cdot)\cdot)\cdot),((\cdot \cdot)(\cdot\cdot))\right) 
Z\left(\triq;((\cdot \cdot)(\cdot\cdot)),(\cdot (\cdot(\cdot\cdot)))\right)$$
$$=$$
$$Z\left(\triq;(((\cdot \cdot)\cdot)\cdot),((\cdot (\cdot\cdot))\cdot)\right)
Z\left(\triq;((\cdot (\cdot\cdot))\cdot),(\cdot ((\cdot\cdot)\cdot))\right)
Z\left(\triq;(\cdot ((\cdot\cdot)\cdot)),(\cdot (\cdot(\cdot\cdot)))\right)$$
The degree k part of the difference between the pentagon equations of $\Phi_Z$ and $\Phi_{Z^K}$ is equivalent to the equation $d(\psi) =0$.
\item[C3] $$\psi^{321}=-\psi$$
Indeed, since $r_v(\assnumi)=\left(\ass\right)^{-1}$ and because of the functoriality and of the invariance under $r_v$, $\Phi_Z$ satisfies 
$$\Phi_Z^{321}=\Phi_Z^{-1}$$
Again, the degree k part of the difference between this equation satisfied by $\Phi_Z$
and the similar one for $\Phi_{Z^K}$ is equivalent to Relation C3.
\item[C2] $$\psi - \psi^{132} -\psi^{213}=0$$
Here, we use the {\em hexagon relation\/} in $\CP_3$
which is obtained by evaluating the good monoidal functor $Z$ on the following tangle
$$ \begin{pspicture}[.2](0,0)(.8,.5)
\psline{-}(.65,0)(0.05,.5)
\psline[border=1pt]{-}(0.05,0)(.45,.5)
\psline[border=1pt]{-}(0.25,0)(.65,.5)
\end{pspicture}
=
\left(2 \times \begin{pspicture}[.2](0,0)(.5,.4)
\psline{-}(0.05,0)(.45,.4)
\rput(.15,-0.05){\tiny 1}
\end{pspicture} \right)_{\ast} 
\left(
\begin{pspicture}[.2](0,0)(.5,.4)
\psline{-}(.45,0)(0.05,.4)
\psline[border=1pt]{-}(0.05,0)(.45,.4)
\rput(.15,-0.05){\tiny 1}
\end{pspicture}
\right)
= \left( \ass \right)
\left( 1_{\itriu} \otimes \fcp \right)
\left(\begin{pspicture}[.2](0,0)(.5,.4)
\psline{-}(0.05,0)(.05,.4)
\rput(.05,-.1){\tiny 1}
\psline{-}(.15,0)(.35,.4)
\rput(.15,-.1){\tiny 3}
\psline{-}(.45,0)(0.45,.4)
\rput(.45,-.1){\tiny 2}
\end{pspicture}
\right)^{-1}
\left(\fcp \otimes 1_{\itriu}\right)
\left( \begin{pspicture}[.2](0,0)(.5,.4)
\psline{-}(0.05,0)(.05,.4)
\rput(.05,-.1){\tiny 3}
\psline{-}(.15,0)(.35,.4)
\rput(.15,-.1){\tiny 1}
\psline{-}(.45,0)(0.45,.4)
\rput(.45,-.1){\tiny 2}
\end{pspicture}
\right)
$$
in the two ways suggested by its two expressions.
The degree k part of the difference between the hexagon equations of $\Phi_Z$ and $\Phi_{Z^K}$ gives rise to the equation $\psi - \psi^{132} + \psi^{312}=0$ 
which is equivalent to C2 thanks to C3.
\item[C4] $$\varepsilon_1(\psi)=\varepsilon_2(\psi)=\varepsilon_3(\psi)=0$$
This is an easy consequence of the fact that $\varepsilon_i(\Phi_Z)=\varepsilon_i(\Phi_{Z^K})=1$.
\end{description}
Now, Proposition 3 of Section 8 in \cite{lm} asserts that for any degree k element $\psi$ of $\CP_3$ which satisfies
C1, C2, C3 and C4, there exists a symmetric degree k element $f$ of $\CP_2$ such that 
$\varepsilon_1(f)=0$ and $d(f)=\psi$. Thus, this proposition applies.
\eop

Now, we can state the following more general uniqueness statement which groups our result and the Le and Murakami theorem.

\begin{theorem}
For any good  monoidal functor $Z$, $a^Z$ is a two-leg element of $\CA(I)$ and there exists a twist $F$ such that $$Z(L)=\Psi(2a^Z)\circ(Z^K)^F(L) $$ for any q-tangle $L$ whose components have as many minima and maxima. For these $L$,  $Z(L)$ can also be written as
$$Z(L)=\left(\Psi(2a^Z)\circ Z^K\right)^{\Psi(2a^Z)(F)}(L).$$ 
where $\Psi(2a^Z)(F)$ is a twist.
\end{theorem}

\section{On the denominators of the Kontsevich integral}
\label{secden}
\setcounter{equation}{0}

Let $M$ be a one-manifold. We say that the {\em denominator\/} of an element $z_n$ of $\CA_n(M)$ {\em divides into\/} an integer $N$
if $Nz_n$ may be written as an integral combination of chord diagrams.  
Of course, the denominator is the greatest common divisor of the $N$ that the denominator divides into. 
Observe that any diagram on $M$ is an integral combination of chord diagrams modulo STU.

When $\Gamma$ is a degree $n$ diagram, $u^{\Gamma}$, $t^{\Gamma}$ and $e^{\Gamma}$ denote the number of univalent vertices, the number of trivalent vertices, and the number of edges of $\Gamma$, respectively. They are related by the obvious equalities $t^{\Gamma}+u^{\Gamma}=2n$ and $2e^{\Gamma}=3t^{\Gamma}+u^{\Gamma}$ which imply:
$$e^{\Gamma} + u^{\Gamma} =3n.$$

This section is devoted to the proof of Corollary~\ref{corcoco} that relies on the following results.

We consider a link $L$ whose components have both zero framing and zero Gauss 
integral. Let  $Z^{CS}$ denote
the perturbative expression of the Chern-Simons theory. Recall that the anomaly vanishes in even degrees, $\alpha_{2i}=0$ for any integer i.

\begin{theorem}[Poirier]
\label{proppoi}
Let $n$ be an integer. Assume $n \geq 3$.
\begin{enumerate}
\item The denominator of the degree n part $Z_n^{CS}(L)$ of 
$Z^{CS}(L)$ divides into $(3n-4)!2^{3n-4}$.
\item The denominator of the degree n part $\alpha_n$ of $\alpha$ divides into $(3n-4)!2^{3n-4}$.
\item $\alpha_n$ belongs to the lattice generated by the
$$ \frac{1}{e^{\Gamma}!2^{e^{\Gamma}}} [\Gamma] $$
where $\Gamma$ runs among the connected degree n diagrams with at least 4 univalent vertices.
\item $\alpha_3=0$.
\end{enumerate}
\end{theorem}
\bp
It follows from \cite[Prop. 1.11, 1.9 and 1.10]{p3} that the denominator
of $Z_n^{CS}(L)$ divides into $(3n-3)!2^{3n-4}$.
In fact, by \cite[Remark~1.12]{p3} whose proof has been explained to the author by Sylvain Poirier and should shortly appear, we can improve this
denominator into $(3n-4)!2^{3n-4}$.
Similarly, the denominator of $2 \alpha_n$ divides into $(3n-4)!2^{3n-4}$ by
\cite[Def.~6.4 and Prop.~6.2]{p3}, and the symmetry of \cite[Lemma 6.6]{p3} allows us to divide this denominator estimate by 2. This is enough to prove the first two assertions for a reader
who knows Poirier's work. Since the denominator estimates are not stated for the anomaly, we now sketch the proofs of the second and the third assertions. (The proof of the first one is similar, but more complicated.)
Now, $n$ is an odd integer greater than 2. In \cite[Def.~6.4 and Prop.~6.2]{p3}, S.~Poirier expresses $2 \alpha_n$ 
as a {\em degree\/} of a map $\Psi$ from a {\em glued configuration space\/} to $\left( S^2 \right)^{3n-3}$. Here, the glued configuration space is an algebraic combination
of smooth oriented manifolds $C_{\Gamma}$ with corners indexed by {\em labelled\/} connected diagrams $\Gamma$ that are degree n diagrams endowed with an orientation of each of their edges and an injection of the set of their edges into  $\{1,2, \dots, 3n-3\}$ (that numbers the edges). The coefficient of a $C_{\Gamma}$ in the combination is 
$$\beta(\Gamma)=\frac{\left(3n-3 -e^{\Gamma}\right)!}{(3n-3)!2^{e^{\Gamma}}} [\Gamma].$$
The map $\Psi$ smoothly maps $C_{\Gamma}$ to $\left( S^2 \right)^{3n-3}$ so that
the regular values of (all the) $\Psi$ form a dense open set and the number of preimages of $\Psi$ (counted with signs and coefficients) of a regular value is constant on this dense set and is, by definition, the degree of $\Psi$. $2 \alpha_n$ is this degree which therefore belongs to the lattice
generated by the $\beta(\Gamma)$ such that $\Psi(C_{\Gamma})$ has interior points. By \cite[Lemma 1.9]{p3}, the involved $\Gamma$ must have at least
4 univalent vertices and therefore at most $(3n-4)$ edges.

Furthermore, when a labelled diagram appears in the above combination, the underlying unlabelled diagram appears with all its possible labellings.

The symmetry of \cite[Lemma 6.6]{p3} shows us that the algebraic preimage of a regular
value of $\Psi$ in some $C_{\Gamma}$ is the same as in $C_{-\Gamma}$ where
$-\Gamma$ is obtained from $\Gamma$ by reversing the orientations of all the edges.

The most delicate point is the content of \cite[Remark~1.12]{p3} that together
with the above arguments allows us to go from $(3n-3)!$ to $(3n-4)!$.
Its proof amounts to choosing a limit generic point of $\left( S^2 \right)^{3n-3}$
with the property that its algebraic preimage in a $C_{\Gamma}$ only depends on the order of the edges induced by the numbering injection and not on the genuine
numbering injection. The existence of such a limit point proves that the degree $2 \alpha_n$ is an integral combination of terms of the form
$$ \frac{2}{e^{\Gamma}!2^{e^{\Gamma}}} [\Gamma] $$
for connected degree $n$ diagrams $\Gamma$ that satisfy $e^{\Gamma} \leq 3n-4$.
The $2$ in the numerator comes from the above symmetry.

The last property $\alpha_3=0$ is \cite[Proposition 1.5]{p3}.
\eop

\begin{lemma}
\label{lempbwo}
Let $k \in \NN$. Let $\pi_k$ denote the composition of the inverse $\chi^{-1}:\CA(I) \longrightarrow \CB$ of the Poincar\'e-Birkhoff-Witt isomorphism by the natural projection from $\CB$ to $\CB_{\ast,k}$. Let $\Gamma$ be a diagram with support $I$ with $u$ univalent vertices.
If $u \geq k$, then $k!(k+1)!(k+2)!\dots(u-1)!\pi_k([\Gamma])$ is an integral combination of k-leg diagrams. If $u<k$, then $\pi_k([\Gamma])=0$.
\end{lemma}
\bp
We assume that the univalent vertices of $\Gamma$ are given the orientation induced by the global orientation of $I$ as in Remark~\ref{rkor}.
Let $\CO(\Gamma)$ be the element of $\CB_{\ast,u}$ obtained from $\Gamma$ by
removing $I$.
We first prove that $$(u-1)!\left([\Gamma] -\chi(\CO(\Gamma))\right)$$
is an integral combination of (classes of) diagrams with $(u-1)$ univalent vertices. 

We use
the orientation of $I$ to number the univalent vertices of $\Gamma$ from $1$ to $u$. When $\sigma$ is a permutation of $\Sigma_u$, $\sigma.\Gamma$ denotes the diagram obtained from
$\Gamma$ by changing the injection from $U$ to $I$ so that the vertex $j$
becomes the $\sigma(j)^{\rm th}$ one on $I$.
Let $\rho$ be the cyclic permutation of $\Sigma_u$ that maps $j$ to $j+1$ mod $u$ for any $i$. It has already been noticed that 
$$[(\rho \sigma). \Gamma]= [\sigma. \Gamma]$$ in $\CA(I)$.

Then $$\chi(\CO(\Gamma))=\frac1{u!}\sum_{\sigma \in \Sigma_u}[\sigma.\Gamma]
=\frac1{(u-1)!}\sum_{\overline{\sigma} \in <\rho>\backslash\Sigma_u}[\sigma.\Gamma]$$
In particular, 
$$[\Gamma]=\chi(\CO(\Gamma)) + \frac1{(u-1)!}\sum_{\overline{\sigma} \in <\rho>\backslash\Sigma_u}\left([\Gamma] - [\sigma.\Gamma]\right)$$
By $STU$, when $\tau$ is a transposition of two consecutive elements $i$ and $i+1$, $\left([\sigma .\Gamma] - [\tau \sigma.\Gamma]\right)$ is equivalent to a diagram
with $(u-1)$ univalent vertices. Now, any permutation $\sigma$ can be written as a product
$\tau_1 \tau_2 \dots \tau_r$ of such transpositions. Therefore
$$[\Gamma] - [\sigma.\Gamma]=\sum_{i=1}^{r}\left([\tau_{i+1} \tau_{i+2} \dots \tau_r.\Gamma]-[\tau_i \tau_{i+1} \tau_{i+2}\dots \tau_r.\Gamma] \right)$$
is an integral combination of diagrams with $(u-1)$ univalent vertices
and $(u-1)!\left([\Gamma] -\chi(\CO(\Gamma))\right)$
is an integral combination of (classes of) diagrams with $(u-1)$ univalent vertices.

Thus, by induction, we have
$$[\Gamma]=\chi(\CO(\Gamma)) + \frac1{(u-1)!}\chi(g_{u-1}) + 
\frac1{(u-2)!(u-1)!}\chi(g_{u-2}) + \dots + \frac1{1!2!\dots(u-2)!(u-1)!} \chi(g_1)$$
where $g_i$ is an integral combination of i-leg diagrams in $\CB_{\ast,i}$.

(In particular, we have proved that $\chi$ is onto. In fact, this proof is built on an excerpt of the Bar-Natan proof of the 
Poincar\'e-Birkhoff-Witt theorem \cite[Theorem 8]{bn}.)
Now, since $\chi$ is an isomorphism, 
$\pi_u([\Gamma])=\CO(\Gamma)$, $\pi_{k}([\Gamma])=\frac1{k!(k+1)!\dots(u-1)!} g_k$
if $k<u$ and $\pi_k([\Gamma])=0$ if $k>u$.
\eop

We say that the {\em two-leg denominator\/} of a two-leg element $z_n$ of $\CA_n(S^1)$ {\em divides into\/} an integer $N$
if $Nz_n$ may be written as an integral combination of two-leg diagrams. 

\begin{lemma}
\label{lemdent}
When $n \geq 5$, the two-leg denominator of $2\alpha_n$ divides into 
$$D(2;n)=(2!3! \dots (n-1)!)\frac{(n-1)!}{(n-4)!}(3n-4)!2^{2n-3}$$
\end{lemma}
\bp
According to Theorem~\ref{proppoi}, $\alpha_n$ belongs to the lattice generated by the
$$ \frac{1}{e^{\Gamma}!2^{e^{\Gamma}}} [\Gamma] $$
where $\Gamma$ runs among the connected degree n diagrams with at least 4 univalent vertices.
For a connected diagram $\Gamma$, $e^{\Gamma} \geq 2n-1$. Since $e^{\Gamma} + u^{\Gamma} =3n$, $u^{\Gamma} \leq n+1$. 
In particular, because of Lemma~\ref{lempbwo}, the 2-leg denominator of $\alpha_n=\chi(\pi_2(\alpha_n))$ divides into the lowest common multiple of the $$d_u=2!3!4!\dots(u-1)!(3n-u)!2^{3n-u}$$
where $u$ belongs to $\{4, 5, \dots,n+1\}$, and it is enough to prove that
$d_u$ divides into $2D(2;n)$
for these $u$.
When $u \leq n$, it is enough to see that $2^{n-u+2}$ divides into $(u!(u+1)! \dots (n-1)!)\frac{(n-1)!}{(n-4)!}$ which is a multiple of $4!^{n-u} \times 4 \times 2$
and thus of $2^{3n-3u+3}$.
When $u=n+1$, 
$$d_{n+1}=2!3!4!\dots(n-1)!n!(2n-1)!2^{2n-1}=(2!3!4!\dots(n-1)!)\frac{(n-1)!}{(n-4)!}(2n)!(n-4)!2^{2n-2}$$ that divides into $2D(2;n)$.
\eop

We will also use the following lemma.

\begin{lemma}[Vogel] 
\label{lemsymvo}
Two-leg elements of $\CB_{\ast,2}$ are symmetric with respect to the exchange of their
two legs.
\end{lemma}
\bp Since a chord is obviously symmetric, we can 
restrict ourselves to 
a two-leg diagram with at least one trivalent vertex and whose two univalent vertices are respectively numbered by 1 and 2. We draw it as 
$$ 
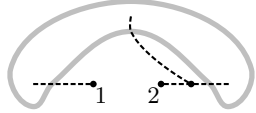
\begin{pspicture}[.2](0,0)(3.4,1.5)
\psccurve[linewidth=2pt,linecolor=lightgray](1.7,1)(2.75,.3)(3,0)(3.3,.5)(1.7,1.4)(.1,.5)(.4,0)
(.65,.3)
\psline[linestyle=dashed,dash=2pt 1pt]{-*}(.4,.3)(1.2,.3)
\rput(1.3,.15){\small 1}
\psline[linestyle=dashed,dash=2pt 1pt]{-*}(3,.3)(2.1,.3)
\rput(2,.15){\small 2}
\pscurve[linestyle=dashed,dash=2pt 1pt]{-*}(1.7,1.2)(1.7,1)(2.5,.3)
 \end{pspicture}$$
where the dashed trivalent part inside the thick topological circle
is not represented. Applying Lemma~\ref{lemcom} where the annulus is a neighborhood of the thick topological circle that contains the pictured trivalent vertex shows that
this diagram is equivalent to 
$$ \begin{pspicture}[.2](0,0)(3.4,1.5)
\psccurve[linewidth=2pt,linecolor=lightgray](1.7,1)(2.75,.3)(3,0)(3.3,.5)(1.7,1.4)(.1,.5)(.4,0)
(.65,.3) 
\psline[linestyle=dashed,dash=2pt 1pt]{-*}(.4,.3)(1.4,.3)
\rput(1.5,.15){\small 1}
\psline[linestyle=dashed,dash=2pt 1pt]{-*}(3,.3)(2.3,.3)
\rput(2.2,.15){\small 2}
\pscurve[linestyle=dashed,dash=2pt 1pt]{-*}(1.7,1.2)(1.7,1)(1,.3)
\end{pspicture}$$
This yields the relations
$$  \begin{pspicture}[.2](.3,0)(2.5,1)
\psline[linestyle=dashed,dash=2pt 1pt]{-*}(.4,.3)(1,.3)
\rput(1.1,.15){\small 1}
\psline[linestyle=dashed,dash=2pt 1pt]{-*}(2.4,.3)(1.6,.3)
\rput(1.5,.15){\small 2}
\pscurve[linestyle=dashed,dash=2pt 1pt]{-*}(1.4,1)(1.4,.8)(2,.3)
 \end{pspicture} =
\begin{pspicture}[.2](.3,0)(2.5,1)
\psline[linestyle=dashed,dash=2pt 1pt]{-*}(.4,.3)(1.2,.3)
\rput(1.3,.15){\small 1}
\psline[linestyle=dashed,dash=2pt 1pt]{-*}(2.4,.3)(1.8,.3)
\rput(1.7,.15){\small 2}
\pscurve[linestyle=dashed,dash=2pt 1pt]{-*}(1.4,1)(1.4,.8)(.8,.3)
 \end{pspicture} =
\begin{pspicture}[.2](.3,0)(2.5,1)
\psline[linestyle=dashed,dash=2pt 1pt]{-*}(1.4,1)(1.4,.6)
\rput(1.55,.7){\small 1}
\psline[linestyle=dashed,dash=2pt 1pt]{-*}(2.4,.3)(1.6,.3)
\rput(1.4,.3){\small 2}
\pscurve[linestyle=dashed,dash=2pt 1pt]{-*}(.4,.3)(1.9,.1)(2,.3)
 \end{pspicture}
=
\begin{pspicture}[.2](.3,0)(2.5,1)
\psline[linestyle=dashed,dash=2pt 1pt]{-*}(.4,.3)(.9,.3)
\rput(1,.15){\small 2}
\psline[linestyle=dashed,dash=2pt 1pt]{-*}(1.4,1)(1.4,.3)
\rput(1.5,.15){\small 1}
\psline[linestyle=dashed,dash=2pt 1pt]{-*}(2.4,.3)(1.4,.6)
 \end{pspicture}
=
\begin{pspicture}[.2](.3,0)(2.5,1)
\psline[linestyle=dashed,dash=2pt 1pt]{-*}(.4,.3)(1,.3)
\rput(1.1,.15){\small 2}
\psline[linestyle=dashed,dash=2pt 1pt]{-*}(2.4,.3)(1.6,.3)
\rput(1.5,.15){\small 1}
\pscurve[linestyle=dashed,dash=2pt 1pt]{-*}(1.4,1)(1.4,.8)(2,.3)
 \end{pspicture}
$$
\eop

Recall that $Z_1^{K}(L)$ is an integer (this is true as soon as $L$ is framed by even integers) and that the denominator of $Z_2^{K}(L)$ is known to be $24$.
Our theorem implies that $$Z^{K}(L)=\Psi(2 \alpha)^{-1}(Z^{CS}(L))$$
where $2 \alpha_1 = \tata$,  $\alpha_3=0$, and $\alpha_{2i}=0$ for any integer i.

In particular, for $i \leq 4$,
$$Z_i^{CS}(L)=Z_i^{K}(L)$$

Now, we want to see how  $\Psi(2 \alpha)^{-1}$ transforms the estimates on the denominators of $Z^{CS}(L)$. 
So, we first make this inverse more explicit. 
To do it we equip $\CB_{\ast,2}$ with the commutative product that
maps a pair $(\Gamma_1,\Gamma_2)$ of two two-leg diagrams to the diagram obtained by identifying one leg of $\Gamma_1$ to one leg of $\Gamma_2$.  We also shift the graduation of $\CB_{\ast,2}$ by -1 so that this product becomes a graded product. Now, \begin{pspicture}[0.2](0,0)(.5,.3)
\psline[linestyle=dashed,dash=2pt 1pt]{*-*}(0.05,.15)(.45,.15)
\end{pspicture} is the degree 0 unit.

When $b$ is an element of $\CB_{\ast,2}$, 
$\firstleg b \secondleg$ denotes the element of $\CA(\{1,2\} \times I)$ obtained by putting one leg of $b$ on $\{1\} \times I$
and the other one on $\{2\} \times I$.
$$2 \firstleg b \secondleg = (2 \times I)_{\ast}(\chi(b)) - \chi(b)_1 -\chi(b)_2 \in \CA(\{1,2\} \times I) $$
Consequently, the denominator of $2 \firstleg b \secondleg$ divides into the denominator
of $\chi(b)$. Obviously, the denominator of $\firstleg b \secondleg$ divides into the two-leg denominator
of $b$.

\begin{lemma}
Let $$A_{2k}=\chi^{-1}(2\alpha_{2k+1})$$
and $$A=\sum_{k \in \NN}A_{2k}$$
Let $B=\sum_{k \in \NN}B_{2k}$ be the element of the algebra $\CB_{\ast,2}$ defined by induction on $k$ by the formula
$$\sum_{k \in \NN}A^{2k+1}B_{2k}=1$$
Then $B$ satisfies the following properties:
\begin{enumerate}
\item $\Psi(\chi(B))$ is the inverse of $\Psi(2 \alpha)$.\\
\item $B_2=0.$\\
\item Let $$D(k)=\left\{\begin{array}{ll}
(6k-1)!2^{6k-1} & \mbox{if}\;1 \leq k \leq 3\\
(2!3! \dots (2k-4)!)(2k-4)!3^2(6k-1)!2^{4k+3} &  \mbox{if}\;k \geq 4
\end{array} \right. $$
the denominator of $\firstleg B_{2k}\secondleg$ divides into $D(k)$ for $k \geq 1$.\\
\item For any positive integers $k_1$ and $k_2$, $D(k_1)D(k_2)$ divides into $D(k_1+k_2)$.
\end{enumerate}
\end{lemma}
\bp
First note that $B_0= A_0=1$ and that $B_{2n}$ is well-defined by induction on $n$ by the given formula in degree $2n$ that reads:
$$B_{2n}+(2n-1)B_{2n-2}A_2+\left((2n-3)A_4 +\left(\begin{array}{c} 2n-3 \\ {2} \end{array}\right)A_2^2\right))B_{2n-4} +\dots+A_{2n}=0$$
and yields $B_2=0$ since $A_2=0$, $B_4=-A_4$ and $B_6=-A_6$.
Applying $\Psi(2 \alpha) \circ \Psi(\chi(B))$ to a chord diagram amounts to inserting $B$ on each chord and then replacing each inserted $B$ by $\sum_{k \in \NN}A^{2k+1}B_{2k}=1$. Thus, $\Psi(2 \alpha) \circ \Psi(\chi(B))=\mbox{Id}$. Since $\Psi(2 \alpha)$ is known to be an isomorphism,
$\Psi(\chi(B))$ is its inverse.

Now, we prove that  the denominator of $\firstleg B_{2k}\secondleg$ divides into $D(k)$ for $k \geq 1$.
 The above defining induction formula for $B_{2k}$, shows that $B_{2k}$ is a homogeneous polynomial of degree $2k$ in the $A_{2j}$ $0<j\leq k$
with integral coefficients, by induction. 

In particular, it is enough to prove that the denominators of the 
\begin{center}
$\firstleg$ degree $2k$ monomials in the $A_{2j}$ $\secondleg$
\end{center}
divide into $D(k)$.
Among these monomials, we have $A_{2k}$, and the denominator of $\firstleg A_{2k} \secondleg$ divides into  $(6k-1)!2^{6k-1}$ (because the denominator of
$2\firstleg A_{2k} \secondleg$
divides into  the denominator of $2 \alpha_{2k+1}$ that divides into the (even) $(6k-1)!2^{6k-2}$ by Theorem~\ref{proppoi}).
When $1 \leq k \leq 3$, we have no other monomial. Thus, the denominator of $\firstleg B_{2k}\secondleg$ divides into  $D(k)=(6k-1)!2^{6k-1}$.

Now, let $k \geq 4$. In order to prove the third assertion, it is enough to prove:
\begin{enumerate}
\item $(6k-1)!2^{6k-1}$ divides into $D(k)$.
\item If $k_2 \geq k_1 \geq 2$, $D(2;2k_1+1)D(2;2k_2+1)$ divides into $D(k_1+k_2)$.
\item If $k_1 \geq 2$ and if $k \geq 4$, $D(2;2k_1+1)D(k)$ divides into $D(k+k_1)$.
\end{enumerate}
where $D(2;.)$ is the two-leg denominator defined in Lemma~\ref{lemdent}:
$$D(2;2k+1)=(2!3!
 \dots (2k)!)\frac{(2k)!}{(2k-3)!}(6k-1)!2^{4k-1}$$
The first and the third assertions are easy to check. Let us check the second
one. Set $k=k_1+k_2$. It is enough to check that
$$(2!3! \dots (2k_1)!)\frac{(2k_1)!}{(2k_1-3)!}(6k_1-1)!(6k_2)!$$
divides into
$$(2k_2+1)!(2k_2+2)!\dots(2k-4)!\frac{(2k-4)!}{(2k_2)!}(2k_2-3)!6k_2(6k-1)!3^22^5.$$
If $k_1\geq 4$, then it is easy to see that 
$$\frac{(5!6! \dots (2k_1)!)}{(2k_1-3)!}(2k_1)!2!3!4! \;\;\mbox{divides into}\;\;
(2k_2+1)!(2k_2+2)!\dots(2k-4)!$$
and we are done.\\
If $k_1=3$, and $k_2 \geq 4$, then we see that 
$2!4!5!6!6!$ divides into $(2k_2+1)!(2k_2+2)!\frac{(2k_2+2)!}{(2k_2)!}(2k_2-3)!$.\\
If $k_1=k_2=3$, then $2!4!4!6!6!$ divides into $6(2k_2+1)!(2k_2+2)!\frac{(2k_2+2)!}{(2k_2)!}(2k_2-3)!$, and 
$5$ divides into $\frac{(6k-1)!}{(6k_2)!(6k_1-1)!}$.\\
If $k_1=2$, since $2!3!4!=3^22^5$, 
it suffices to check that $4!$ divides into 
$(2k_2-3)!6k_2 \frac{(6k-1)!}{(6k_2)!(6k_1-1)!}$.
Thus, we have proved that  the denominator of $\firstleg B_{2k}\secondleg$ divides into $D(k)$, and the last assertion of the lemma that relies on the same arguments is easy to check.
\eop

\noindent{\sc Proof of Corollary~\ref{corcoco}:}
 $$Z^{K}(L)=\Psi(\chi(B))(Z^{CS}(L))$$
Let $n \geq 4$. First assume that all the components of $L$ are framed by zero. We want to prove that the denominator of $Z^{K}_n(L)$ divides into $$d(n)=\left\{\begin{array}{ll}
(3n-4)!2^{3n-4} & {\rm if}\;3 \leq n \leq 8\\
(2!3! \dots (n-5)!)(n-5)!3^2(3n-4)!2^{2n+1} &  {\rm if}\;n \geq 9 \;\mbox{and n is odd}\\
(2!3! \dots (n-6)!)(n-6)!3^2(3n-4)!2^{2n-1} &  {\rm if}\;n \geq 10 \;\mbox{and n is even}
\end{array} \right. $$
It is enough to prove that for any $r\leq n$ and for any set $\{2k_1, \dots, 2k_j\}$ of even integers $2k_i$ greater than 3, such that $ 2\sum_{i=1}^j k_i= n-r$, and $j \leq r$, the product
$$\mbox{denominator of } Z_r^{CS}(L)\prod_{i=1}^j\mbox{denominator 
of } \firstleg B_{2k_i} \secondleg$$ divides into $d(n)$.
Thus, since $\prod_{i=1}^jD(k_i)$ divides into $D\left(\sum_{i=1}^jk_i\right)$, it is enough to prove
that $(3n-4)!2^{3n-4}$ divides into $d(n)$,
and that for any even integer $2k$ greater than $3$,
\begin{enumerate}
\item $D(k)$ divides into $d(2k+1)$,
\item $24 D(k)$ divides into $d(2k+2)$, and,
\item For any $r \geq 3$, $(3r-4)!2^{3r-4}D(k)$ divides into $d(2k+r)$.
\end{enumerate}
Recall that
$$D(k)=\left\{\begin{array}{ll}
(6k-1)!2^{6k-1} & {\rm if}\;1 \leq k \leq 3\\
(2!3! \dots (2k-4)!)(2k-4)!3^2(6k-1)!2^{4k+3} &  {\rm if}\;k \geq 4
\end{array} \right. $$
that makes the checking easy.
A framing change on a component induces a multiplication of $Z^{K}(L)$ by some
$\exp(\frac{k}{2}\tata)$, with $k \in \ZZ$, on this component that does not change the denominator
estimate.
\eop

\begin{remark}
The above denominators could be improved by using the result $\alpha_5=0$, but 
the author has not checked the unpublished Poirier proof of this result that involves a Maple computation.
\end{remark}

\end{document}